\documentclass[MSNbibl,nameyear,dvips]{arxstspdf}
\usepackage{mathbh}
\usepackage{graphicx}
\usepackage{flushend}
\usepackage{stfloats}

\volume{28}
\issue{3}
\pubyear{2013}
\firstpage{289}
\lastpage{312}
\doi{10.1214/13-STS434} 

\makeatletter

\renewcommand{\citep}[1]{(\citeauthor{#1}, \citeyear{#1})}
\newtheorem{theorem}{Theorem}

\newtheorem{lemma}[theorem]{Lemma}

\newproclaim{definition}[theorem]{Definition}
\newproclaim{eg}[theorem]{Example}
\newproclaim{example}[theorem]{Example}

\newcommand{\bern}{\operatorname{Bern}}
\newcommand{\tb}{\operatorname{Beta}}
\newcommand{\disc}{\operatorname{Discrete}}
\newcommand{\ga}{\operatorname{Ga}}
\newcommand{\pois}{\operatorname{Pois}}
\newcommand{\ppp}{\operatorname{PPP}}
\newcommand{\unif}{\operatorname{Unif}}

\newcommand{\mbe}{\mathbb{E}}
\newcommand{\mbp}{\mathbb{P}}
\newcommand{\mbo}{\mathbh{1}}

\newcommand{\eqd}{\stackrel{d}{=}}
\newcommand{\iid}{\stackrel{\mathrm{i.i.d}.}{\sim}}
\newcommand{\indep}{\stackrel{\mathrm{indep}}{\sim}}

\newcommand{\rmloc}{\phi}
\newcommand{\rmsp}{\Phi}
\newcommand{\rmlocdist}{H}
\newcommand{\rmweight}{\xi}
\newcommand{\randm}{\mu}
\newcommand{\intens}{\nu}
\newcommand{\basedistr}{\rmlocdist}
\newcommand{\likedistr}{\mathcal{L}}
\newcommand{\data}{X}
\newcommand{\detpart}{\pi} 
\newcommand{\randpart}{\Pi} 
\newcommand{\blockpart}{A} 
\newcommand{\randlabpart}{\rmloc} 
\newcommand{\eppf}{p} 
\newcommand{\indpart}{Z} 
\newcommand{\partblockfreq}{\rho} 
\newcommand{\detfa}{f} 
\newcommand{\randfa}{F} 
\newcommand{\randofa}{\tilde{F}} 
\newcommand{\detofa}{\tilde{f}} 
\newcommand{\ok}{\tilde{K}} 
\newcommand{\randfasp}{\mathcal{F}} 
\newcommand{\blockfa}{\blockpart} 
\newcommand{\randlabfa}{\randlabpart} 
\newcommand{\restfa}{\mathcal{R}} 
\newcommand{\efpf}{\eppf} 
\newcommand{\indfa}{Y} 
\newcommand{\fablockfreq}{\rmweight} 
\newcommand{\perm}{\sigma} 
\newcommand{\stickprop}{V} 
\newcommand{\omsprop}{W} 
\newcommand{\soms}{w} 
\newcommand{\urnextra}{\kappa} 
\newcommand{\subord}{T} 
\newcommand{\invsub}{S} 
\newcommand{\drift}{c} 
\newcommand{\levym}{\Lambda} 
\newcommand{\levydens}{\rho} 
\newcommand{\sumdens}{f} 
\newcommand{\subjumpsize}{\rmweight} 
\newcommand{\subjumploc}{\rmloc} 
\newcommand{\lapexp}{\Psi} 
\newcommand{\sumofjumps}{\tau} 
\newcommand{\dpconc}{\theta} 
\newcommand{\dproc}{D} 
\newcommand{\gpconc}{\dpconc} 
\newcommand{\gpscale}{b} 
\newcommand{\gproc}{G} 
\newcommand{\bpmass}{\gamma} 
\newcommand{\bpconc}{\dpconc} 
\newcommand{\indhasf}{I} 
\newcommand{\bproc}{B} 

\makeatother

\begin{document}
\begin{frontmatter}

\title{Cluster and Feature Modeling from Combinatorial Stochastic Processes}
\runtitle{Cluster and Feature Modeling}

\begin{aug}
\begin{aug}
\author[a]{\fnms{Tamara} \snm{Broderick}\corref{}\ead[label=e1]{tab@stat.berkeley.edu}},
\author[b]{\fnms{Michael I.} \snm{Jordan}\ead[label=e2]{jordan@stat.berkeley.edu}}
\and
\author[c]{\fnms{Jim} \snm{Pitman}\ead[label=e3]{pitman@stat.berkeley.edu}}
\end{aug}
\runauthor{T. Broderick, M.~I. Jordan and J. Pitman}

\affiliation{University of California, Berkeley}

\address[a]{T. Broderick is Graduate Student, Department of Statistics, University of California,
Berkeley, Berkeley, California 94720, USA \printead{e1}.}
\address[b]{M.~I. Jordan is Pehong Chen Distinguished Professor, Department of EECS and Department of
Statistics, University of California, Berkeley, Berkeley, California 94720, USA
\printead{e2}.}
\address[c]{J. Pitman is Professor, Department of Statistics and Department of Mathematics, University of
California, Berkeley, Berkeley, California 94720, USA \printead{e3}.}
\end{aug}

%
\begin{abstract}
One of the focal points of the modern literature on Bayesian nonparametrics
has been the problem of \emph{clustering}, or \emph{partitioning},
where each data point is modeled as being associated with one and only
one of some collection of groups called clusters or partition blocks.
Underlying these Bayesian nonparametric models are a set of
interrelated stochastic processes, most notably the Dirichlet process
and the Chinese restaurant process. In this paper we provide a formal
development of an analogous problem, called \emph{feature modeling},
for associating data points with arbitrary nonnegative integer numbers
of groups, now called features or topics. We review the existing
combinatorial stochastic process representations for the clustering
problem and develop analogous representations for the feature modeling
problem. These representations include the beta process and the Indian
buffet process as well as new representations that provide insight into
the connections between these processes. We thereby bring the same
level of completeness to the treatment of Bayesian nonparametric
feature modeling that has previously been achieved for Bayesian
nonparametric clustering.
\end{abstract}

%
\begin{keyword}
\kwd{Cluster}
\kwd{feature}
\kwd{Dirichlet process}
\kwd{beta process}
\kwd{Chinese restaurant process}
\kwd{Indian buffet process}
\kwd{nonparametric}
\kwd{Bayesian}
\kwd{combinatorial stochastic process}\vspace*{-5pt}
\end{keyword}

\end{frontmatter}

\section{Introduction}

Bayesian nonparametrics is the area of Bayesian analysis in which
the finite-dimensional prior distributions of classical Bayesian
analysis are replaced with stochastic processes. While the rationale
for allowing infinite collections of random variables into Bayesian
inference is often taken to be that of diminishing the role of prior
assumptions, it is also possible to view the move to nonparametrics
as supplying the Bayesian paradigm with a richer collection of
distributions with which to express prior belief, thus in some sense
emphasizing the role of the prior. In practice, however, the field has
been dominated by two stochastic processes---the Gaussian process and the
Dirichlet process---and thus the flexibility promised by the nonparametric
approach has arguably not yet been delivered. In the current paper
we aim to provide a broader perspective on the kinds of stochastic
processes that can provide a useful toolbox for Bayesian nonparametric
analysis. Specifically, we focus on \emph{combinatorial stochastic
processes} as embodying mathematical structure that is useful for both
model specification and inference.

The phrase ``combinatorial stochastic process''\break comes from probability
theory \citep{pitman:2006:combinatorial}, where it refers to
connections between
stochastic processes and the mathematical field of combinatorics.
Indeed, the focus in this area of probability theory is on random
versions of classical combinatorial objects such as partitions, trees
and graphs---and on the role of combinatorial analysis in establishing
properties of these processes. As we wish to argue, this connection is
also fruitful in a statistical setting. Roughly speaking, in statistics
it is often natural to model observed data as arising from a combination
of underlying factors. In the Bayesian setting, such models are often
embodied as latent variable models in which the latent variable has
a compositional structure. Making explicit use of ideas from
combinatorics in latent variable modeling cannot only suggest
new modeling ideas but can also provide essential help with
calculations of marginal and conditional probability distributions.

The Dirichlet process already serves as one interesting exhibit of the
connections
between Bayesian nonparametrics and combinatorial stochastic processes.
On the one hand, the Dirichlet process is classically defined in terms of
a partition of a probability space \citep{Ferguson73}, and there are many
well-known connections between the Dirichlet process and urn
models (\cite{blackwell:1973:ferguson}; \cite{hoppe:1984:polya}).
In the current paper, we will review and expand upon some of these
connections, beginning our treatment (nontraditionally) with the
notion of an \emph{exchangeable partition probability function} (EPPF)
and, from there, \mbox{discussing} related urn models, stick-breaking representations,
subordinators and random measures.

On the other hand, the Dirichlet process is limited in terms of the
statistical notion of a ``combination of underlying factors'' that we
referred to above. Indeed, the Dirichlet process is generally used
in a statistical setting to express the idea that each data point is
associated with one and only one underlying factor. In contrast to such
\emph{clustering models}, we wish to also \mbox{consider} \emph{featural models},
where each data point is associated with a set of underlying features
and it is the interaction among these features that gives rise to
an observed data point. Focusing on the case in which these features
are binary, we develop some of the combinatorial stochastic process
machinery needed to specify featural priors. Specifically, we develop
a counterpart to the EPPF,\vadjust{\goodbreak} which we refer to as the \emph{exchangeable
feature probability function} (EFPF), that characterizes the combinatorial
structure of certain featural models. We again develop connections between
this combinatorial function and suite of related stochastic processes,
including urn models, stick-breaking representations, subordinators
and random measures. As we will discuss, a particular underlying random measure
in this case is the \emph{beta process}, originally studied by \citet
{hjort:1990:nonparametric}
as a model of random hazard functions in survival analysis, but adapted
by \citet{ThibauxJo07} for applications in featural modeling.

For statistical applications it is not enough to develop expressive
prior specifications, but it is also essential that inferential
computations involving the posterior distribution are tractable.
One of the reasons for the popularity of the Dirichlet process is
that the associated urn models and stick-breaking representations
yield a variety of useful inference algorithms \citep{neal:2000:markov}.
As we will see, analogous algorithms are available for featural models.
Thus, as we discuss each of the various representations associated
with both the Dirichlet process and the beta process, we will also
(briefly) discuss some of the consequences of each
for posterior inference.

The remainder of the paper is organized as follows.
We start by reviewing partitions and introducing feature
allocations in Section~\ref{sec:partition_feature} in order to
define distributions over these models (Section~\ref{sec:epf}) via the
EPPF
in the partition case (Section~\ref{sec:epf_eppf}) and the
EFPF in the
feature allocation case (Section~\ref{sec:epf_efpf}).
Illustrating these exchangeable
probability functions with examples, we will see that the well-known
\textit{Chinese restaurant process} (CRP) \citep{aldous:1985:exchangeability}
corresponds to a particular EPPF choice (Example~\ref{ex:epf_crp})
and the \textit{Indian buffet process} (IBP) (Griffiths and Ghahramani, \citeyear{griffiths:2006:infinite})
corresponds to a particular choice of EFPF (Example~\ref{ex:epf_ibp}).
From here, we progressively build up richer models by first reviewing
stick lengths (Section~\ref{sec:stick}), which we will see represent limiting
frequencies of certain clusters or features, and then subordinators
(Section~\ref{sec:sub}), which further associate a random label with
each cluster
or feature. We illustrate these progressive augmentations for both the CRP
(Examples~\ref{ex:epf_crp},~\ref{ex:cond_crp},~\ref{ex:stick_crp},
\ref{ex:sub_crp} and~\ref{ex:sub_crp_sticks})
and IBP examples (Examples~\ref{ex:epf_ibp},~\ref{ex:cond_ibp}, \ref
{ex:stick_ibp} and~\ref{ex:sub_ibp}).
We augment the model once more to obtain a random measure on a general space
of cluster or feature parameters in Section~\ref{sec:crm}, and discuss how
marginalization
of this random measure yields the CRP\vadjust{\goodbreak} in the case of the Dirichlet process
(Example~\ref{ex:crm_dp}) and the IBP in the case of the beta process
(Example~\ref{ex:crm_bp}).
Finally, in Section~\ref{sec:conclusion}, we mention some of the other
combinatorial
stochastic processes, beyond the Dirichlet process and the beta process,
that have begun to be studied in the Bayes\-ian nonparametrics literature,
and we provide suggestions for further developments.

\section{Partitions and Feature Allocations} \label{sec:partition_feature}

While we have some intuitive ideas about what constitutes
a cluster or feature model, we want to formalize these ideas before
proceeding. We begin with the underlying combinatorial structure
on the data indices. We think of $[N] := \{1,\ldots,N\}$ as
representing the
indices of the first $N$ data points. There are different groupings
that we apply in the cluster case (\textit{partitions}) and feature case
(\textit{feature allocations}); we describe these below.

First, we wish to describe the space of \textit{partitions} over the
indices $[N]$.
In particular, a partition $\detpart_{N}$ of $[N]$ is defined to be a
collection of
mutually exclusive, exhaustive, nonempty subsets of $[N]$ called
\textit{blocks}; that is,
$\detpart_{N} = \{\blockpart_{1},\ldots,\blockpart_{K}\}$ for some
number of partition blocks $K$.
An example partition of $[6]$ is $\detpart_{6} = \{\{1,3,4\}, \{2\}, \{
5,6\}\}$.
Similarly, a partition of $\mathbb{N} = \{1,2,\ldots\}$ is a
collection of
mutually exclusive, exhaustive, nonempty subsets of $\mathbb{N}$.
In this case, the number of blocks may be infinite, and we have
$\detpart_{N} = \{A_{1},A_{2},\ldots\}$. An example partition of
$\mathbb{N}$ into
two blocks is
$\{\{n\dvtx n \mbox{ is even}\}, \{n\dvtx n \mbox{ is odd}\}\}$.

We introduce a generalization of a partition called a \textit{feature
allocation} that relaxes
both the mutually exclusive and exhaustive restrictions. In particular,
a feature allocation $\detfa_{N}$ of $[N]$ is defined to be a multiset
of nonempty subsets of $[N]$, again called \textit{blocks},
such that each index $n$ can belong to any finite number of blocks.
Note that the
constraint that no index should belong to infinitely many blocks
coincides with our intuition for the meaning of these blocks
as groups to which the index belongs.
Consider an example where the data points are images
expressed
as pixel arrays, and the latent features represent animals that
may or may not appear in each
picture. It is impossible to display an infinite number of animals in
a picture with finitely many pixels.

We write $\detfa_{N} = \{\blockfa_{1},\ldots,\blockfa_{K}\}$ for some
number of feature
allocation blocks $K$.
An example feature allocation of $[6]$ is
$\detfa_{6} = \{\{2,3\},\{2,4,6\},\{3\},\{3\},\{3\}\}$. Just as the
blocks of a partition
are sometimes called \textit{clusters}, so are the blocks of a feature
allocation
sometimes called \textit{features}.
We note that a partition is always a feature allocation, but the
converse statement
does not hold in general; for instance, $f_{6}$ given above is not a partition.

In the remainder of this section we continue our development in terms
of feature allocations since partitions are a special case of the
former object.
We note that we can extend the idea of
random partitions \citep{aldous:1985:exchangeability}
to consider \textit{random feature allocations}. If $\randfasp_{N}$
is the space of all feature allocations of $[N]$, then a random feature
allocation $\randfa_{N}$ of $[N]$
is a random element of this space.

We next introduce a few useful assumptions on our random
feature allocation.
Just as exchangeability of observations is often a central assumption in
statistical modeling, so will we make use of \textit{exchangeable feature
allocations}. To rigorously define such feature allocations, we
introduce the
following notation. Let $\perm\dvtx  \mathbb{N} \rightarrow\mathbb{N}$
be a
finite permutation. That is, for some finite value $N_{\perm}$, we have
$\perm(n) = n$ for all
$n > N_{\sigma}$. Further, for any block $\blockfa\subset\mathbb{N}$,
denote the permutation applied to the block as follows:
$\perm(\blockfa) := \{\perm(n)\dvtx n \in\blockfa\}$. For any feature
allocation $F_{N}$,
denote the permutation applied to the feature allocation as follows:
$\perm(\randfa_{N}) := \{\perm(\blockfa)\dvtx \blockfa\in\randfa
_{N}\}$.
Finally, let $\randfa_{N}$ be a random feature allocation of $[N]$.
Then we say that
$\randfa_{N}$ is exchangeable if $\randfa_{N} \eqd\sigma(\randfa_{N})$
for every finite permutation $\perm$.

Our second assumption in what follows will be that we are dealing with a
\textit{consistent} feature allocation. We often implicitly imagine the
indices arriving one at a time:
first 1, then 2, up to $N$ or beyond. We will find it useful,
similarly, in
defining random feature allocations to suppose that the randomness at
stage $n$ somehow agrees with the randomness at stage $n+1$. More formally,
we say that a feature allocation $\detfa_{M}$ of $[M]$ is a \emph{restriction} of
a feature allocation $\detfa_{N}$ of $[N]$ for $M < N$ if
\[
\detfa_{M} = \bigl\{\blockfa\cap[M]\dvtx \blockfa\in
\detfa_{N}\bigr\}.
\]
Let $\restfa_{N}(\detfa_{M})$ be the set of all feature allocations
of $[N]$
whose restriction to $[M]$ is $\detfa_{M}$. Then we say that the sequence
of random feature allocations $(\randfa_{n})$ is \textit{consistent} if
for all $M$ and
$N$ such that $M < N$, we have that
%
%
\begin{equation}
\label{eq:strong_consistency} \randfa_{N} \in
\restfa_{N}(\randfa_{M})\quad \mbox{a.s.}
\end{equation}

With this consistency condition in hand, we can define a random feature
allocation
$\randfa_{\infty}$
of $\mathbb{N}$. In particular, such a feature allocation is
characterized by the sequence of
consistent finite restrictions $\randfa_{N}$ to $[N]$:
$\randfa_{N} := \{A \cap[N]\dvtx A \in\randfa_{\infty}\}$. Then
$\randfa
_{\infty}$
is equivalent to a consistent sequence of finite feature allocations
and may be thought
of as a random
element of the space of such sequences: $\randfa_{\infty} = (\randfa
_{n})_{n}$.
We let $\randfasp_{\infty}$ denote the space of consistent feature
allocations, of which
each random feature allocation is a random element, and we see that the
sigma-algebra
associated with this space is generated by the finite-dimensional sigma-algebras
of the restricted random feature allocations $\randfa_{n}$.

We say that $\randfa_{\infty}$ is exchangeable if
$\randfa_{\infty} \eqd\sigma(\randfa_{\infty})$ for every finite
permutation
$\perm$. That is, when the permutation $\perm$ changes no indices above
$N$, we
require $\randfa_{N} \eqd\sigma(\randfa_{N})$, where $\randfa_{N}$
is the
restriction of $\randfa_{\infty}$ to $[N]$. A characterization of
distributions
for $\randfa_{\infty}$ is provided by
\citet{broderick:2013:feature}, where a similar treatment of the
introductory ideas
of this section also appears.

%
\begin{figure*}

\includegraphics{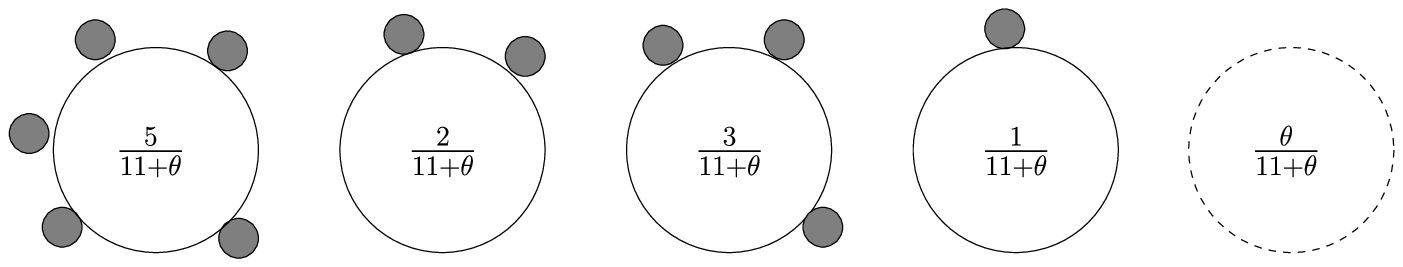}

\caption{The diagram represents a possible CRP seating
arrangement after 11 customers have entered a restaurant with parameter
$\dpconc$. Each large white circle is a table, and the smaller gray
circles are customers sitting at those tables. If a 12th customer
enters, the expressions in the middle of each table give the
probability of the new customer sitting there. In particular, the
probability of the 12th customer sitting at the first table is
$5/(11+\dpconc)$, and the probability of the 12th customer forming a
new table is $\dpconc/(11+\dpconc)$.}\label{fig:crp}
\end{figure*}

In what follows, we consider particular useful ways
of representing distributions for exchangeable, consistent random
feature allocations with emphasis on partitions as a special case.

\section{Exchangeable Probability Functions} \label{sec:epf}

Once we know that we can construct (exchangeable
and consistent) random partitions and feature allocations,
it remains to find useful representations of distributions over these
objects.

\subsection{Exchangeable Partition Probability Function} \label{sec:epf_eppf}

Consider first an exchangeable, consistent, random partition
$(\randpart_{n})$.
By the exchangeability
assumption, the distribution of the partition should depend only on the
(unordered) sizes of the
blocks. Therefore, there exists a function
$\eppf$ that is symmetric in its arguments such that, for any specific
partition assignment
$\detpart_{n} = \{\blockpart_{1},\ldots,\blockpart_{K}\}$, we have
%
%
\begin{equation}
\label{eq:eppf} \mbp(\randpart_{n} = \detpart_{n}) =
\eppf\bigl(|\blockpart_{1}|, \ldots, |\blockpart_{K}|\bigr).
\end{equation}
The function $\eppf$ is called the \textit{exchangeable partition
probability function}
(EPPF) \citep{pitman:1995:exchangeable}.

%
%
\begin{example}[(Chinese restaurant process)] \label{ex:epf_crp}
The Chinese restaurant process (CRP) (Blackwell and MacQueen, \citeyear{blackwell:1973:ferguson})
is an iterative description of a partition
via the conditional distributions of the partition blocks to which increasing
data indices belong.
The Chinese restaurant\vadjust{\goodbreak}
metaphor forms an equivalence between customers entering a Chinese restaurant
and data indices; customers who share a table at the restaurant represent
indices belonging to the same partition block.

To generate the label for the first index,
the first customer enters the restaurant and sits
down at some table, necessarily unoccupied since no one else is in the
restaurant. A ``dish'' is set out at the new table; call the dish ``1''
since it is
the first dish.
The customer is assigned the label of the dish at her table: $\indpart
_{1} = 1$.
Recursively, for a restaurant with \textit{concentration parameter}
$\theta$,
the $n$th customer sits at an occupied table with probability in proportion
to the number of people at the table and at a new table with
probability proportional to $\dpconc$.
In the former case, $\indpart_{n}$ takes the
value of the existing dish at the table, and, in the latter case, the next
available dish $k$ (equal to the number of existing tables plus one) appears
at the new table, and $\indpart_{n} = k$. By summing over all
possibilities when the
$n$th customer arrives, one obtains the normalizing constant for the
distribution across potential occupied tables: $(n-1+\dpconc)^{-1}$.
An example of the distribution over tables for the $n$th customer
is shown in Figure~\ref{fig:crp}.
To summarize,
if we let $K_{n} := \max\{\indpart_{1},\ldots,\indpart_{n}\}$, then
the distribution of table assignments for the $n$th customer is
%
%
\begin{eqnarray}\label{eq:crp}
&& \mbp(\indpart_{n} = k | \indpart_{1},
\ldots,\indpart_{n-1})
\nonumber
\\[-8pt]
\\[-8pt]
\nonumber
 &&\quad= (n-1+\dpconc)^{-1}
\cases{\#\{m\dvtx m < n, \indpart_{m} = j\},\vspace*{2pt}\cr
 \quad\hspace*{11pt} \mbox{for }j \le
K_{n-1},
\vspace*{2pt}\cr
\theta,\quad\mbox{for }k = K_{n-1}+1.}
\end{eqnarray}

We note that an equivalent generative description follows a
P\'olya urn style in specifying that each incoming customer
sits next to an existing customer with probability proportional to 1 and
forms a new table with probability proportional to
$\dpconc$ \citep{hoppe:1984:polya}.

Next, we find the probability of the partition induced by considering
the collection of indices sitting at each table as a block in the
partition. Suppose that $N_{k}$ individuals sit at table $k$ so that
the set of cardinalities of nonzero table occupancies is $\{
N_{1},\ldots
,N_{K}\}$ with $N := \sum_{k=1}^{K} N_{k}$. That is, we are considering
the case when $N$ customers have entered the restaurant and sat at $K$
different tables in the specified configuration.

We can see from equation (\ref{eq:crp}) that when the $n$th customer
enters ($n > 1$),
we obtain a factor of $n - 1 + \dpconc$ in the denominator. Using the
following notation for the rising and falling factorial
\[
x_{M \uparrow a} := \prod_{m=0}^{M-1} (x +
ma),\quad x_{M \downarrow
a} := \prod_{m=0}^{M-1} (x
- ma),
\]
we find a factor of $(\dpconc+ 1)_{N - 1 \uparrow1}$ must occur in
the denominator of the probability of the partition of~$[N]$.
Similarly, each time a customer forms a new table except for the first
table, we obtain a factor of $\dpconc$ in the numerator. Combining
these factors, we find a factor of $\dpconc^{K-1}$ in the numerator.
Finally, each time a customer sits at an existing table with $n$
occupants, we obtain a factor of $n$ in the numerator. Thus, for each
table $k$, we have a factor of $(N_{k} - 1)!$ once all customers have
entered the restaurant.\looseness=1

Having collected all terms in the process, we see that the probability
of the resulting configuration is
%
%
\begin{equation}
\label{eq:crp_eppf} \mbp(\randpart_{N} =
\detpart_{N}) = \frac{\dpconc^{K-1} \prod_{k=1}^{K} (N_{k} - 1)!}{ (\dpconc+
1)_{N-1 \uparrow1} }.
\end{equation}
We first note that equation (\ref{eq:crp_eppf}) depends only on the
block sizes and
not on the order of arrival of the customers or dishes at the tables.
We conclude that the partition generated according to the CRP scheme is
exchangeable. Moreover, as the partition $\randpart_{M}$ is the
restriction of $\randpart_{N}$ to $[M]$ for any $N > M$ by
construction, we have that equation (\ref{eq:crp_eppf}) satisfies the
consistency
condition. It follows that equation (\ref{eq:crp_eppf}) is, in fact,
an EPPF.
\end{example}
%
%

\subsection{Exchangeable Feature Probability Function} \label{sec:epf_efpf}

Just as we considered an exchangeable, consistent, random partition
above, so we now
turn to an exchangeable, consistent, random feature allocation
$(\randfa_{n})$.
Let $\detfa_{N} = \{\blockpart_{1},\ldots,\blockpart_{K}\}$ be any
particular feature
allocation.
In calculating $\mbp(\randfa_{N} = \detfa_{N})$, we start by demonstrating
in the next example that
this probability in some sense undercounts features when they contain exactly
the same indices: for example, $A_{j} = A_{k}$ for some $j \ne k$.
For instance, consider the following example.

%
%
\begin{example}[(A two-block, Bernoulli feature allocation)] \label{ex:two_bern}
Let $q_{A}, q_{B} \in(0,1)$ represent the frequencies of features $A$
and $B$.
Draw $\indpart_{A,n} \iid\bern(q_{A})$ and $\indpart_{B,n} \iid
\bern
(q_{B})$, independently.
Construct the random feature allocation by collecting those
indices with successful draws:
\[
\randfa_{N} := \bigl\{\{n\dvtx n \le N, \indpart_{A,n} = 1
\}, \{n\dvtx n \le N, \indpart_{B,n} = 1\}\bigr\}.
\]
Then the probability of the feature allocation $\randfa_{5} = \detfa
_{5} := \{\{2,3\},\{2,3\}\}$ is
\[
q_{A}^{2} (1-q_{A})^{3}
q_{B}^{2} (1-q_{B})^{3},
\]
but the probability of the feature allocation $\randfa_{5} = \detfa
'_{5} := \{\{2,3\},\{2,5\}\}$ is
\[
2 q_{A}^{2} (1-q_{A})^{3}
q_{B}^{2} (1-q_{B})^{3}.
\]
The difference is that in the latter case the features can be
distinguished, and so
we must account for the two possible pairings of features to
frequencies $\{q_{A},q_{B}\}$.

Now, instead, let $\randofa_{N}$ be $\randfa_{N}$ with a uniform
random ordering
on the features. There is just a single possible ordering of $\detfa
_{5}$, so
the probability of $\randofa_{5} = \detofa_{5} := (\{2,3\}, \{2,3\})$
is again
\[
q_{A}^{2} (1-q_{A})^{3}
q_{B}^{2} (1-q_{B})^{3}.
\]
However, there are two orderings of $\detfa'_{5}$, so the
probability of $\randofa_{5} = \detofa'_{5} := (\{2,5\}, \{2,3\})$ is
\[
q_{A}^{2} (1-q_{A})^{3}
q_{B}^{2} (1-q_{B})^{3},
\]
and the same holds for the other ordering.
\end{example}

For reasons suggested by the previous example, we will find it useful
to work with
the random feature allocation after uniform random
ordering, $\randofa_{N}$. One way to achieve such an ordering
and maintain consistency across different $N$ is to associate some
independent, continuous
random variable with each feature; for example, assign a uniform random
variable on $[0,1]$ to
each feature and order the features according to the order of
the assigned random variables.
When we view feature allocations constructed as marginals of
a \textit{subordinator} in Section~\ref{sec:sub}, we will see that
this construction
is natural.

In general,
given a probability of a random feature allocation, $\mbp(\randfa_{N} =
\detfa_{N})$,
we can find the probability of a \textit{random ordered feature allocation},
$\mbp(\randofa_{N} = \detofa_{N})$ as follows. Let $H$ be the
number of unique elements of $\randfa_{N}$, and let $(\ok_{1},\ldots
,\ok_{H})$
be the multiplicities of these unique elements in decreasing size. Then
%
%
\begin{equation}
\label{eq:order_mult} \mbp(\randofa_{N} =
\detofa_{N}) = \pmatrix{ K
\cr
\ok_{1}, \ldots,
\ok_{H} }^{-1} \mbp(\randfa _{N} =
\detfa_{N}),\hspace*{-25pt}
\end{equation}
where
\[
\pmatrix{ K
\cr
\ok_{1}, \ldots, \ok_{H} } :=
\frac{ K! }{ \ok
_{1}! \cdots
\ok_{H}!}.
\]

We will see in Section~\ref{sec:sub} that augmentation of an exchangeable
partition with a random ordering is also natural.
However, the probability of an ordered random partition
is not substantively different from the probability of an unordered
version since the
factor contributed by ordering a partition is always $1 / K!$, where $K$
here is the number of partition blocks.

With this framework in place, we can see that
some ordered feature allocations have a probability function $\efpf$
nearly as in equation (\ref{eq:eppf}),
that is, moreover, symmetric in its block-size arguments. Consider
again the previous example.

%
%
\begin{example}[(A two-block, Bernoulli feature allocation (continued))]
Consider any $\randfa_{N}$ with block sizes $N_{1}$ and $N_{2}$ constructed
as in Example~\ref{ex:two_bern}. Then
%
%
\begin{eqnarray}\label{eq:efpf_two_bern}
\nonumber
&&\mbp(\randofa_{N} = \detofa_{N})\\
&&\quad =
\tfrac{1}{2} q_{A}^{N_{1}} (1-q_{A})^{N - N_{1}}
q_{B}^{N_{2}} (1-q_{B})^{N - N_{2}}
\nonumber
\\
&&\qquad {} + \tfrac{1}{2} q_{A}^{N_{2}}
(1-q_{A})^{N - N_{2}} q_{B}^{N_{1}}
(1-q_{B})^{N - N_{1}}
\nonumber\\
 &&\quad= \efpf(N, N_{1},
N_{2}),
\end{eqnarray}
where $\efpf$ is some function of the number of indices $N$
and the block sizes $(N_{1}, N_{2})$ that we
note is symmetric in all arguments after the first. In particular,
we see that the order of $N_{1}$ and $N_{2}$ was immaterial.
\end{example}

We note that in the partition case, $\sum_{k=1}^{K} |\blockpart_{k}| =
N$, so
$N$ is implicitly an argument to the EPPF. In the feature case, this
summation condition no longer holds,
so we make the argument $N$ explicit in equation (\ref{eq:efpf_two_bern}).

However, it is not necessarily the case that such a function, much less
a symmetric one, exists for exchangeable feature models---in contrast to
the case of exchangeable partitions and the EPPF.

%
%
\begin{example}[(A general two-block feature allocation)]
We here describe an exchangeable, consistent random feature allocation
whose (ordered) distribution does not depend only on the number of
indices $N$
and the sizes of the blocks of the feature
allocation.

Let $p_{1}, p_{2}, p_{3}, p_{4}$ be fixed frequencies that sum to one.
Let $\indfa_{n}$ represent the collection of features to which index
$n$ belongs.
For $n \in\{1,2\}$, choose $\indfa_{n}$ independently and identically
according to
\[
\indfa_{n} = \cases{ %
 \{A\} ,&
$\mbox{with probability } p_{1},$
\vspace*{2pt}\cr
\{B\},& $\mbox{with probability } p_{2},$
\vspace*{2pt}\cr
\{A,B\}, & $\mbox{with probability } p_{3},$
\vspace*{2pt}\cr
\varnothing,& $\mbox{with probability } p_{4}.$}
\]
We form a feature allocation from these labels as follows. For each
label ($A$ or $B$),
collect those indices $n$ with the given label appearing in $\indfa
_{n}$ to form a feature.

Now consider two possible outcome feature allocations:
$\detfa_{2} = \{\{2\}, \{2\}\}$ and $\detfa'_{2} = \{\{1\},\{2\}\}$.
The likelihood of any random ordering $\detofa_{2}$ of $\detfa_{2}$
under this model is
\[
\mbp(\randofa_{2} = \detofa_{2}) = p_{1}^{0}
p_{2}^{0} p_{3}^{1}
p_{4}^{1}.
\]
The likelihood of any ordering $\detofa'_{2}$ of $\detfa'_{2}$ is
\[
\mbp\bigl(\randofa_{2} = \detofa'_{2}\bigr) =
p_{1}^{1} p_{2}^{1}
p_{3}^{0} p_{4}^{0}.
\]
It follows from these
two likelihoods that we can choose values of $p_{1},p_{2},p_{3},p_{4}$ such
that $\mbp(\randofa_{2} = \detofa_{2}) \ne\mbp(\randofa_{2} =
\detofa
'_{2})$. But
$\detofa_{2}$ and $\detofa'_{2}$ have the same block counts and $N$
value ($N = 2$). So there can be no such symmetric function $\eppf$,
as in equation (\ref{eq:efpf_two_bern}), for this model.
\end{example}

When a function $\efpf$ exists in the form
%
%
\begin{equation}
\label{eq:efpf} \mbp(\randofa_{N} = \detofa_{N}) =
\efpf\bigl(N,|\blockfa_{1}|, \ldots, |\blockfa_{K}|\bigr)
\end{equation}
for some random ordered feature allocation
$
\detofa_{N} = (A_{1},\ldots,A_{K})
$
such that $\efpf$ is symmetric in all arguments after the first, we
call it the \emph{exchangeable feature probability function} (EFPF). Note that
the EPPF is not a special\vadjust{\goodbreak} case of the EFPF. The EPPF assigns zero
probability to any multiset in which an index occurs in more than one
element of the
multiset; only the sizes of the multiset blocks are relevant in the
EFPF case.

We next consider a more complex example of an EFPF.

%
%
\begin{example}[(Indian buffet process)] \label{ex:epf_ibp}
The Indian buffet process (IBP) (Griffiths and Ghahramani,\break \citeyear{griffiths:2006:infinite})
is a generative model for a random feature allocation that is specified
recursively like the Chinese restaurant process. Also like the CRP,
this culinary metaphor forms an equivalence between customers and the
indices $n$ that will be partitioned: $n \in\mathbb{N}$. Here,
``dishes'' again correspond to feature labels just as they corresponded
to partition labels for the CRP. But in the IBP case, a customer can
sample multiple dishes.

In particular, we start with a single customer, who enters the buffet
and chooses $K^{+}_{1} \sim\pois(\bpmass)$ dishes. Here, $\bpmass> 0$
is called the \textit{mass parameter}, and we will also see the
\emph{concentration parameter} $\bpconc> 0$ below. None of the dishes have
been sampled by any other customers since no other customers have yet
entered the restaurant. We label the dishes $1,\ldots,K^{+}_{1}$ if
$K^{+}_{1} > 0$. Recursively, the $n$th customer chooses which dishes
to sample in two parts. First, for each dish $k$ that has previously
been sampled by any customer in $1,\ldots,n-1$, customer $n$ samples
dish $k$ with probability $N_{n-1,k} / (\bpconc+ n - 1)$ for $N_{n,k}$
equal to the number of customers indexed $1,\ldots,n$ who have tried
dish $k$. As each dish represents a feature, and sampling a dish
represents that the customer index $n$ belongs to that feature,
$N_{n,k}$ is the size of the block of the feature labeled $k$ in the
feature allocation of $[n]$. Next, customer $n$ chooses $K^{+}_{n} \sim
\pois(\bpconc\bpmass/(\bpconc+n-1))$ new dishes to try. If $K^{+}_{n}
> 0$, then the dishes receive unique labels $K_{n-1}+1, \ldots, K_{n}$.
Here, $K_{n}$ represents the number of sampled dishes after $n$
customers: $K_{n} = K_{n-1} + K^{+}_{n}$.
An example of the first few steps in the Indian buffet process is shown
in Figure~\ref{fig:ibp}.

%
\begin{figure}

\includegraphics{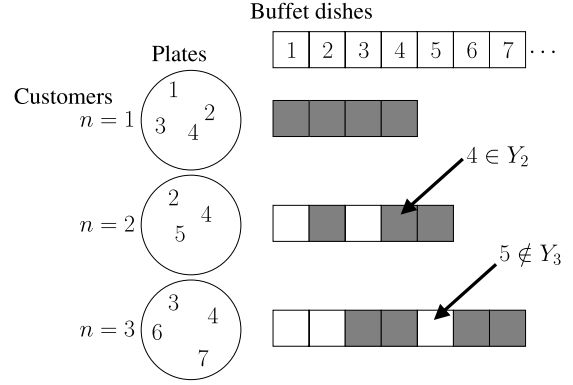}

\caption{Illustration of an Indian buffet process. The
buffet (\textit{top}) consists of a vector of dishes, corresponding to
features. Each customer---corresponding to a data point---who enters
first decides whether or not to eat dishes that the other customers
have already sampled and then tries a random number of new dishes, not
previously sampled by any customer. A gray box in position $(n,k)$
indicates customer $n$ has sampled dish $k$, and a white box indicates
the customer has not sampled the dish. In the example, the second
customer has sampled exactly those dishes indexed by 2, 4 and 5:
$\indfa
_{2} = \{2,4,5\}$.}\label{fig:ibp}\vspace*{5pt}
\end{figure}

With this generative model in hand, we can find the probability of a
particular feature allocation. We discover its form by enumeration as
for the CRP EPPF in Example~\ref{ex:epf_crp}. At each round $n$, we
have a Poisson
number of new features, $K^{+}_{n}$, represented. The probability
factor associated with these choices is a product of Poisson densities:
\[
\prod_{n=1}^{N} \frac{1}{K^{+}_{n}!} \biggl(
\frac{\bpconc\bpmass
}{\bpconc+ n - 1} \biggr)^{K^{+}_{n}} \exp \biggl( -\frac{\bpconc
\bpmass
}{\bpconc+ n -1} \biggr).
\]
Let $M_{k}$ be the round on which the $k$th dish, in order of
appearance, is first chosen. Then the denominators for future dish
choice probabilities are the factors in the product $(\bpconc+ M_{k})
\cdot(\bpconc+ M_{k} + 1) \cdots(\bpconc+ N -\nobreak 1)$. The numerators
for the times when the dish is chosen are the factors in the product $1
\cdot\nobreak2 \cdots\break(N_{N,k}-\nobreak1)$. The numerators for the times when the dish
is not chosen yield $(\bpconc+ M_{k} - 1) \cdots(\bpconc+ N - 1 - N_{N,k})$.
Let $\blockfa_{n,k}$ represent the collection of indices in the feature
with label $k$ after $n$ customers have entered the restaurant. Then
$N_{n,k} = |\blockfa_{n,k}|$.
Finally, let $\ok_{1},\ldots,\ok_{H}$ be the multiplicities of unique
feature blocks formed by this model. We note that there are
\[
\Biggl[ \prod_{n=1}^{N}
K^{+}_{n}! \Biggr] \Bigg/ \Biggl[ \prod
_{h=1}^{H} \ok _{h}! \Biggr]
\]
rearrangements of the features generated by this process that all yield
the same feature allocation. Since they all have the same generating
probability, we simply multiply by this factor to find the feature
allocation probability.
Multiplying all factors together and taking $\detfa_{n} = \{\blockfa
_{N,1}, \ldots, \blockfa_{N,K_{N}}\}$ yields
\begin{eqnarray*}
&&\mbp(\randfa_{N} = \detfa_{N})
\\
&&\quad= \frac{ \prod_{n=1}^{N} K_{n}^{+}! }{ \prod_{h=1}^{H} \ok_{h}! }\\
&&\qquad{} \cdot \Biggl[ \prod_{n=1}^{N}
\frac{1}{K^{+}_{n}!} \biggl(\frac
{\bpconc
\bpmass}{\bpconc+ n - 1} \biggr)^{K^{+}_{n}} \exp \biggl( -
\frac
{\bpconc
\bpmass}{\bpconc+ n -1} \biggr) \Biggr]
\\
&&\qquad{} \cdot \Biggl[ \prod_{k=1}^{K_{N}}
\frac{\Gamma(\bpconc+
M_{k})}{\Gamma(\bpconc+ N)} \Gamma(N_{N,k}) \frac{\Gamma(\bpconc+ N - N_{N,k})}{\Gamma(\bpconc
+M_{k}-1)} \Biggr]
\\
&&\quad= \Biggl( \prod_{h=1}^{H}
\ok_{h}! \Biggr)^{-1} \Biggl[ \prod
_{n=1}^{N} (\bpconc\bpmass )^{K^{+}_{n}} \exp
\biggl( -\frac
{\bpconc
\bpmass}{\bpconc+ n -1} \biggr) \Biggr] \\
&
&\qquad{}\cdot \biggl[ \frac{\prod_{k=1}^{K_{N}} (\theta+ M_{k} - 1)}{\prod_{n=1}^{N} (\theta+ n -
1)^{K_{n}^{+}}}
\biggr]
\\
&&\qquad{} \cdot \Biggl[ \prod_{k=1}^{K_{N}}
\frac{\Gamma(N_{N,k}) \Gamma(\bpconc+ N - N_{N,k})}{\Gamma(\bpconc
+ N)} \Biggr]
\\
&&\quad= \Biggl( \prod_{h=1}^{H}
\ok_{h}! \Biggr)^{-1} (\bpconc\bpmass )^{K_{N}} \\
&&\qquad{}\cdot\exp
\Biggl( -\bpconc\bpmass\sum_{n=1}^{N} (
\bpconc+ n - 1)^{-1} \Biggr)\\
&&\qquad{}\cdot \prod_{k=1}^{K_{N}}
\frac{\Gamma(N_{N,k}) \Gamma
(N-N_{N,k}+\bpconc)}{\Gamma(N+\bpconc)}.
\end{eqnarray*}

It follows from equation (\ref{eq:order_mult}) that the probability of
a uniform
random ordering of the feature allocation is
%
%
\begin{eqnarray}
\label{eq:ibp_efpf} && \mbp(\randofa_{N} =
\detofa_{N})\nonumber
\\
&&\quad= \frac{1}{K_{N}!} (\bpconc
\bpmass)^{K_{N}} \exp \Biggl( -\bpconc \bpmass\sum
_{n=1}^{N} (\bpconc+ n - 1)^{-1} \Biggr)
\\
&&\qquad {} \cdot\prod_{k=1}^{K_{N}}
\frac{\Gamma(N_{N,k}) \Gamma
(N-N_{N,k}+\bpconc)}{\Gamma(N+\bpconc)}.\nonumber
\end{eqnarray}

The distribution of $\randofa_{N}$ has no dependence on the ordering of
the indices in $[N]$. Hence,
the distribution of $\randfa_{N}$ depends only on the same
quantities---the number of indices and the feature block sizes---and
the feature multiplicities. So we see that the IBP construction yields
an exchangeable random feature allocation. Consistency follows from the
recursive construction and exchangeability. Therefore, equation (\ref
{eq:ibp_efpf}) is
seen to be in EFPF form [cf. equation~(\ref{eq:efpf})].
\end{example}

Above, we have seen two examples of how specifying a conditional
distribution for the block membership of index $n$ given the block
membership of indices in $[n-1]$ yields an exchangeable probability
function, for example, the EPPF in the CRP case (Example \ref
{ex:epf_crp}) and the
EFPF in the IBP case (Example~\ref{ex:epf_ibp}). This conditional
distribution is
often called a \textit{prediction rule}, and study
of the prediction rule in the clustering case may be referred to as
\textit{species sampling} (\cite{pitman:1996:some};
\cite{hansen:1998:prediction}; \cite{lee:2008:defining}).
We will see next that the prediction rule can conversely be recovered
from the exchangeable probability function specification and,
therefore, the two are equivalent.

\subsection{Induced Allocations and Block Labeling} \label{sec:labels}

In Examples~\ref{ex:epf_crp} and~\ref{ex:epf_ibp} above, we formed
partitions and
feature allocations in the following way. For partitions, we assigned
labels $\indpart_{n}$ to each index $n$. Then we generated a partition
of $[N]$ from the sequence $(\indpart_{n})_{n=1}^{N}$ by saying that
indices $m$ and $n$ are in the same partition block ($m \sim n$) if and
only if $\indpart_{n} = \indpart_{m}$. The resulting partition is
called the \textit{induced partition} given the labels $(\indpart
_{n})_{n=1}^{N}$. Similarly, given labels $(\indpart
_{n})_{n=1}^{\infty
}$, we can form an induced partition of $\mathbb{N}$. It is easy to
check that, given a sequence $(\indpart_{n})_{n=1}^{\infty}$, the
induced partitions of the subsequences $(\indpart_{n})_{n=1}^{N}$ will
be consistent.

In the feature case, we first assigned label collections $\indfa_{n}$
to each index $n$. $\indfa_{n}$ is interpreted as a set containing the
labels of the features to which $n$ belongs. It must have finite
cardinality by our definition of a feature allocation. In this case, we
generate a feature allocation on $[N]$ from the sequence $(\indfa
_{n})_{n=1}^{N}$ by first letting $\{\rmloc_{k}\}_{k=1}^{K}$ be the set
of unique values in $\bigcup_{n=1}^{N} \indfa_{n}$. Then the features
are the collections of indices with shared labels: $\detfa_{N} =
\{\{n\dvtx \rmloc_{k} \in\indfa_{n}\}\dvtx k = 1,\ldots,K\}$. The resulting feature
allocation $\detfa_{N}$ is called the \textit{induced feature allocation}
given the labels $(\indfa_{n})_{n=1}^{N}$. Similarly, given label
collections $(\indfa_{n})_{n=1}^{\infty}$, where each $\indfa_{n}$ has
finite cardinality, we can form an induced feature allocation of
$\mathbb{N}$. As in the partition case, given a sequence $(\indfa
_{n})_{n=1}^{\infty}$, we can see that the induced feature allocations
of the subsequences $(\indfa_{n})_{n=1}^{N}$ will be consistent.

In reducing to a partition or feature allocation from a set of labels,
we shed the information
concerning the labels for each partition block or feature. Conversely,
we introduce
\textit{order-of-appearance} labeling schemes to give partition blocks or
features labels when we have, respectively, a partition or
feature allocation.

In the partition case, the order-of-appearance labeling scheme assigns
the label 1
to the partition block containing index 1. Recursively, suppose we have
seen $n$
indices in $k$ different blocks with labels $\{1,\ldots,k\}$. And
suppose the $n+1$st
index does not belong to an existing block. Then we assign its block
the label $k+1$.

In the feature allocation case, we note that index 1 belongs to
$K_{1}^{+}$ features.
If $K_{1}^{+} = 0$, there are no features to label yet. If $K_{1}^{+} >
0$, we
assign these $K_{1}^{+}$ features labels in $\{1,\ldots,K_{1}^{+}\}$.
Unless otherwise specified, we
suppose that the labels are chosen uniformly at random.
Let $K_{1} = K_{1}^{+}$.
Recursively, suppose we have seen $n$ indices and $K_{n}$ different features
with labels $\{1,\ldots,K_{n}\}$. Suppose the $n+1$st index belongs to
$K_{n+1}^{+}$ features
that have not yet been labeled. Let $K_{n+1} = K_{n} + K_{n+1}^{+}$. If
$K_{n+1}^{+} = 0$, there are no new features to label.
If $K_{n+1}^{+} > 0$, assign these $K_{n+1}^{+}$ features labels in $\{
K_{n} + 1, \ldots, K_{n+1}\}$, for example, uniformly at random.

We can use these labeling schemes to find the prediction rule,
which makes use of partition block and feature labels, from the EPPF or
EFPF as appropriate.
First, consider a partition with EPPF $\eppf$. Then, given labels
$(\indpart_{n})_{n=1}^{N}$ with
$K_{N} = \max\{\indpart_{1},\ldots,\indpart_{N}\}$,
we wish to find the distribution of the label $\indpart_{N+1}$. Using
an order-of-appearance
labeling, we know that either $\indpart_{N+1} \in\{\indpart
_{1},\ldots
,\indpart_{N}\}$ or
$\indpart_{N+1} = K_{N}+1$. Let $\detpart_{N} = \{\blockpart
_{N,1},\ldots,\blockpart_{N,K_{N}}\}$ be the partition
induced by $(\indpart_{n})_{n=1}^{N}$.
Let $N_{N,k} = |\blockpart_{N,k}|$. Let $\mbo(A)$ be the indicator of
event $A$; that is, $\mbo(A)$ equals 1 if $A$ holds and 0 otherwise.
Let $N_{N+1,k} = N_{k} + \mbo\{\indpart_{N+1} = k\}$ for $k=1,\ldots
,K_{N+1}$, and set $N_{N,K_{N}+1} = 0$ for completeness. $K_{N+1} =
K_{N} + \mbo\{\indpart_{N+1} > K_{N}\}$ is the number of partition
blocks in the partition of $[N+1]$.
Then the conditional distribution satisfies
\begin{eqnarray*}
&&\mbp(\indpart_{N+1} = z | \indpart_{1},\ldots,
\indpart_{N}) \\
&&\quad= \frac{\mbp(\indpart_{1},\ldots,\indpart_{N}, \indpart_{N+1} =
z)}{\mbp(\indpart_{1},\ldots,\indpart_{N})}.
\end{eqnarray*}
But the probability of a certain labeling is just the probability of
the underlying
partition in this construction, so
\begin{eqnarray*}
&&\mbp(\indpart_{N+1} = z | \indpart_{1},\ldots,
\indpart_{N})\\
&&\quad   = \frac{\eppf(N_{N+1,1},\ldots,N_{N+1,K_{N+1}})}{\eppf
(N_{N,1},\ldots
,N_{N,K_{N}})}.
\end{eqnarray*}

%
%
\begin{example}[(Chinese restaurant process)] \label{ex:cond_crp}
We continue our Chinese restaurant process example
by deriving the Chinese restaurant table assignment scheme from the
EPPF in equation (\ref{eq:crp_eppf}). Substituting in the EPPF for the
CRP, we find
%
%
\begin{eqnarray}\label{eq:pred_crp_derived}
& &\hspace*{-4pt}\mbp(Z_{N+1} = z | Z_{1},
\ldots,Z_{N})
\nonumber\\
&&\hspace*{-4pt}\quad= \frac{\eppf(N_{N,1},\ldots,N_{N+1,K_{N+1}})}{\eppf
(N_{N,1},\ldots
,N_{N,K_{N}})}
\nonumber
\\
&&\hspace*{-4pt}\quad=
\Biggl( \dpconc^{K_{N+1}-1} \prod_{k=1}^{K_{N+1}} (N_{N+1,k} - 1)!
\Biggr)
\nonumber
\\
&&\hspace*{-4pt}\qquad{}\cdot  \bigl( (\dpconc+ 1)_{(N+1)-1 \uparrow1}  \bigr)^{-1}\nonumber\\
&&\hspace*{-4pt}\qquad{}\bigg /\Biggl(
\Biggl( \dpconc^{K_{N}-1} \prod_{k=1}^{K_{N}} (N_{N,k} - 1)!
\Biggr)\nonumber
\\
&&\hspace*{62pt}{}\cdot\bigl( (\dpconc+ 1)_{N-1 \uparrow1}  \bigr)^{-1}\Biggr)
\nonumber\\
&&\hspace*{-4pt}\quad= (N +
\dpconc)^{-1} \cases{ %
 N_{N,k}, &
$\mbox{for $z = k \le K_{N}$}$,
\vspace*{2pt}\cr
\theta,& $\mbox{for $z = K_{N}+1$,}$ }
\end{eqnarray}
just as in equation (\ref{eq:crp}).
\end{example}

To find the feature allocation prediction rule, we now imagine a feature
allocation with EFPF $\efpf$. Here we must be slightly more careful
about counting
due to feature multiplicities.
Suppose that after $N$ indices have
been seen, we have label collections $(\indfa_{n})_{n=1}^{N}$,
containing a
total of $K_{N}$ features, labeled $\{1,\ldots,K_{N}\}$. We wish to
find the
distribution of $\indfa_{N+1}$. Suppose $N+1$ belongs to $K_{N+1}^{+}$
features that do not contain any index in $[N]$.
Using an order-of-appearance labeling, we know that, if $K_{N+1}^{+} > 0$,
the $K_{N+1}^{+}$ new features have labels $K_{N}+1,\ldots,K_{N} +
K_{N+1}^{+}$.
Let $\detfa_{N} = \{\blockfa_{1},\ldots,\blockfa_{K_{N}}\}$ be the
feature allocation
induced by $(\indfa_{n})_{n=1}^{N}$. Let $N_{N,k} = |\blockfa_{N,k}|$
be the size
of the $k$th feature.
So $N_{N+1,k} = N_{N,k} + \mbo\{k \in\indfa_{N+1}\}$, where we let
$N_{K_{N}+j} = 0$ for
all of the features that are first exhibited by index $N+1$:
$j \in\{1,\ldots,K_{N+1}^{+}\}$. Further, let the number of features,
including new ones,
be written $K_{N+1} = K_{N} + K_{N+1}^{+}$.
Then the conditional distribution satisfies
\begin{eqnarray*}
\mbp(\indfa_{n+1} = y | \indfa_{1},\ldots,
\indfa_{N}) = \frac{\mbp(\indfa_{1},\ldots,\indfa_{N}, \indfa_{N+1} =
y)}{\mbp
(\indfa_{1},\ldots,\indfa_{N})}.
\end{eqnarray*}
As we assume that the labels $\indfa$ are consistent\break across~$N$, the
probability
of a certain labeling is just the probability of
the underlying ordered feature allocation times a combinatorial term.
The combinatorial term accounts first for the uniform ordering of the new
features among themselves for labeling
and then for the uniform ordering of the new features among the old
features in the overall uniform random ordering:
%
%
\begin{eqnarray}\label{eq:pred_from_efpf}
&&\mbp(\indfa_{N+1} = y | \indfa_{1},\ldots,
\indfa_{N})\nonumber
\\
&&\quad= \frac{1}{K_{N+1}^{+}!} \cdot \bigl[
(K_{N}+1) \cdot(K_{N}+2) \cdots K_{N+1} \bigr]
\nonumber\\
&&\qquad {} \cdot\frac{
\efpf(N, N_{N+1,1},\ldots,N_{N+1,K_{N+1}})
}{
\efpf(N, N_{N,1},\ldots,N_{N,K_{N}})
}
\nonumber\\
 &&\quad= \frac{1}{K_{N+1}^{+}!} \cdot
\frac{K_{N+1}!}{K_{N}!}\nonumber\\
&&\qquad{} \cdot\frac{
\efpf(N, N_{N+1,1},\ldots,N_{N+1,K_{N+1}})
}{
\efpf(N, N_{N,1},\ldots,N_{N,K_{N}})
} .
\end{eqnarray}

%
%
\begin{example}[(Indian buffet process)] \label{ex:cond_ibp}
Just as we derived the Chinese restaurant process prediction rule
[equation (\ref{eq:pred_crp_derived})] from
its EPPF [equation (\ref{eq:crp_eppf})] in Example~\ref{ex:cond_crp},
so can we derive the Indian
buffet process prediction rule
from its EFPF [equation (\ref{eq:ibp_efpf})] by using equation (\ref
{eq:pred_from_efpf}).
Substituting the IBP EFPF into equation (\ref{eq:pred_from_efpf}), we
find
\begin{eqnarray*}
&& \mbp(\indfa_{n+1} = y | \indfa_{1},\ldots,
\indfa_{N})
\\
&&\quad= \frac{1}{K_{N+1}^{+}!} \cdot\frac{K_{N+1}!}{K_{N}!}
\biggl(\frac{1}{K_{N+1}!}\biggr)
(\bpconc\bpmass)^{K_{N+1}}\\
&&\qquad {} \cdot \exp \Biggl( -\bpconc\bpmass\sum_{n=1}^{N+1} (\bpconc+ n - 1)^{-1}  \Biggr)\\
&&\qquad {} \cdot\Biggl[\prod_{k=1}^{K_{N+1}}
{\Gamma(N_{N+1,k}) \Gamma\bigl((N+1)-N_{N+1,k}+\bpconc\bigr)}\\
&&\hspace*{120pt}\qquad {} /{\bigl(\Gamma
\bigl((N+1)+\bpconc\bigr)\bigr)}\Biggr]\\
&&\qquad\Bigg/\Biggl\{\biggl(\frac{1}{K_{N}!}\biggr)
(\bpconc\bpmass)^{K_{N}} \\
&&\hspace*{30pt}\quad{}\cdot\exp \Biggl( -\bpconc\bpmass\sum_{n=1}^{N}
(\bpconc+ n - 1)^{-1}  \Biggr)\\
&&\hspace*{32pt}\quad{}\cdot \Biggl[\prod_{k=1}^{K_{N}} {\Gamma
(N_{N,k}) \Gamma(N-N_{N,k}+\bpconc)}\\
&&\hspace*{112pt}\qquad{}/{\bigl(\Gamma(N+\bpconc)\bigr)}\Biggr]\Biggr\}
\\
&&\quad= \biggl[ \frac{1}{K_{N+1}^{+}!} \exp \biggl(- \frac{\bpconc\bpmass}{ \theta+ (N+1) - 1}
\biggr)
\\
&&\hspace*{57pt}{}\cdot\biggl(\frac{\bpconc\bpmass}{ \theta+ (N+1) - 1} \biggr)^{K_{N+1}^{+}} \biggr]
\\
&&\qquad {} \cdot\bigl(\bpconc+ (N+1) - 1\bigr)^{K_{N+1}^{+}}\\
&&\qquad{} \cdot \Biggl[ \prod
_{k=K_{N}+1}^{K_{N+1}} \bigl(\bpconc+ (N+1) - 1
\bigr)^{-1} \Biggr]
\\
&&\qquad {} \cdot\prod_{k=1}^{K_{N}}
\frac{N_{k}^{\mbo\{k \in z\}} (N -
N_{N,k} + \bpconc)^{\mbo\{k \notin z\}} }{N+\bpconc}
\\
&&\quad= \pois \biggl( K_{N+1}^{+} \Big| \frac{\bpconc\bpmass}{ \theta+ (N+1) -
1} \biggr)\\
&&\qquad{}
\cdot\prod_{k=1}^{K_{N}} \bern \biggl( \mbo\{k
\in z\}\Big | \frac
{N_{N,k}}{N + \bpconc} \biggr).
\end{eqnarray*}
The final line is exactly the Poisson distribution for the number of
new features times the Bernoulli distributions for the draws of
existing features, as described in Example~\ref{ex:epf_ibp}.
\end{example}
%
%

\subsection{Inference} \label{sec:epf_infer}

The prediction rule formulation of the EPPF or EFPF is particularly
useful in providing a means of inferring partitions and feature allocations
from a data set. In particular, we assume that we have data points
$\data_{1},\ldots,\data_{N}$ generated in the following manner. In the
partition case,
we generate an exchangeable, consistent, random partition $\randpart_{N}$
according to the distribution specified by some EPPF $\eppf$. Next, we
assign each partition block a random parameter that characterizes that block.
To be precise, for the $k$th partition block to appear according to
an order-of-appearance labeling scheme, give this block a new \textit{random}
label $\randlabpart_{k} \sim\basedistr$, for some continuous distribution
$\basedistr$. For each $n$, let $\indpart_{n} = \randlabpart_{k}$ where
$k$ is the order-of-appearance label of index $n$. Finally, let
%
%
\begin{equation}
\label{eq:likelihood} \data_{n} \indep\likedistr(
\indpart_{n})
\end{equation}
for some distribution $\likedistr$ with parameter
$\indpart_{n}$. The choices of both $\basedistr$ and $\likedistr$
are specific
to the problem domain.

Without attempting to survey the vast literature on clustering,
we describe a stylized example to provide intuition for the preceding
generative model. In this example, let $n$ index an animal observed
in the wild; $\indpart_{n} = \indpart_{m}$ indicates that animals $n$
and $m$ belong to the same (latent, unobserved) species;
$\indpart_{n} = \indpart_{m} = \randlabpart_{k}$ is a vector describing
the (latent, unobserved) height and weight for that species; and
$\data_{n}$ is the observed height and weight of the $n$th animal.

$\data_{n}$ need not even be directly observed, but equation~(\ref
{eq:likelihood}) together
with an EPPF might be part of a larger generative model. In a generalization
of the previous stylized example, $\indpart_{n}$ indicates the dominant
species in the $n$th geographical region; $\indpart_{n} = \randlabpart_{k}$
indicates some overall species height and weight parameters (for the
$k$th species); $\data_{n}$
indicates the height and weight parameters for species $k$ in the $n$th
region. That is, the height and weight for the species may vary by region.
We measure and observe the height and weight $(E_{n,j})_{j=1}^{J}$
of some $J$ animals in the $n$th region, believed to be i.i.d. draws from a
distribution depending on $\data_{n}$.

Note that the sequence $(\indpart_{n})_{n=1}^{N}$ is sufficient to
describe the
partition $\Pi_{N}$ since $\Pi_{N}$ is the collection of blocks of
$[N]$ with the same
label values $\indpart_{n}$. The continuity of $\basedistr$ is
necessary to
guarantee the a.s. uniqueness of the block values.
So, if we can describe the posterior distribution of
$(\indpart_{n})_{n=1}^{N}$,
we can in principle describe the posterior distribution of $\Pi_{N}$.

The posterior distribution of $(\indpart_{n})_{n=1}^{N}$
conditional on $(\data_{n})_{n=1}^{N}$ cannot typically be solved for in
closed form, so we turn to a method that approximates this posterior.
We will see that prediction rules facilitate the design of a Markov
Chain Monte
Carlo (MCMC) sampler, in which we approximate the desired posterior
distribution by a Markov chain of random samples proven to
have the true posterior as its equilibrium distribution.

In the Gibbs sampler formulation of MCMC (Geman and Geman, \citeyear{geman:1984:stochastic}),
we sample each
parameter in turn and conditional on all other parameters in the model.
In our case, we will sequentially sample each element of
$(\indpart_{n})_{n=1}^{N}$. The key observation here is that
$(\indpart_{n})_{n=1}^{N}$ is an exchangeable sequence.
This observation follows by noting that the partition is exchangeable
by assumption, and the sequence
$(\randlabpart_{k})$ is exchangeable since it is i.i.d.;
$(\indpart_{n})$ is an exchangeable sequence since it is a
function of $(\randpart_{n})$ and $(\randlabpart_{k})$.
Therefore, the distribution of $\indpart_{n}$, given the remaining
elements $\mathbf{\indpart}_{-n} :=
(\indpart_{1},\ldots,\indpart_{n-1},\indpart_{n+1},\ldots,\indpart
_{N})$, is the same as
if we thought of $\indpart_{n}$ as the final, $N$th element in a sequence
with $N-1$ preceding values given by $\mathbf{\indpart}_{-n}$.
And the distribution of $\indpart_{N}$ given $\mathbf{\indpart}_{-N}$
is provided by the prediction rule.
The full details of the
Gibbs sampler for the CRP in Examples~\ref{ex:epf_crp} and~\ref{ex:cond_crp}
were introduced by
\citet{escobar:1994:estimating},
\citet{maceachern:1994:estimating},
\citet{escobar:1995:bayesian}
and are covered in fuller generality by \citet{neal:2000:markov}.

It is worth noting that
the sequence of order-of-appearance labels is not exchangeable;
for instance, the first label is always 1.\vadjust{\goodbreak} However, the prediction rule
for $\indpart_{N}$ given $(\indpart_{1},\ldots,\indpart_{N-1})$
breaks into two parts:
(1) the probability of $\indpart_{N}$ taking either a value
in $\{\indpart_{1},\ldots,\indpart_{N-1}\}$ or a new value and
(2) the distribution of $\indpart_{N}$ when it takes a new value.
When programming such a sampler, it is often useful to simply
encode the sets of unique values, which may be done by
retaining any set of labels that induce the correct partition (e.g.,
integer labels) and separately retaining the
set of unique parameter values. Indeed, updating the parameter
values and partition block assignments separately can lead to
improved mixing of the sampler \citep{maceachern:1994:estimating}.

Similarly, in the feature case,
we imagine the following generative model for our data.
First, let $\randfa_{N}$ be a random feature allocation
generated according to the EFPF $\efpf$. For the $k$th
feature block in an order-of-appearance labeling scheme,
assign a random label $\randlabfa_{k} \sim\basedistr$
to this block for some continuous distribution~$\basedistr$.
For each $n$, let $\indfa_{n} = \{\randlabfa_{k}\dvtx k \in J_{n}\}$,
where $J_{n}$ is here the set of order-of-appearance
labels of the features to which $n$ belongs. Finally,
as above,
\[
\data_{n} \indep\likedistr(\indfa_{n}),
\]
where the likelihood $\likedistr$ and parameter
distribution $\basedistr$ are again application-specific
and where now $\likedistr$ depends on the variable-size
collection of parameters in~$\indfa_{n}$.

\citet{griffiths:2011:indian} provide a review of likelihoods used
in practice for feature models. To motivate some of these modeling
choices, let us consider some stylized examples that provide helpful
intuition. For example, let $n$ index customers at a book-selling website;
$\randlabfa_{k}$ describes a book topic such as economics, modern art
or science fiction. If $\randlabfa_{k}$ describes science fiction books,
$\randlabfa_{k} \in\indfa_{n}$ indicates that the $n$th customer
likes to
buy science fiction books. But $\indfa_{n}$ might have cardinality greater
than one (the customer is interested in multiple book topics) or cardinality
zero (the customer never buys books). Finally, $X_{n}$ is a set of
book sales for customer $n$ on the book-selling site.

As a second example, let $n$ index pictures in a database;
$\randlabfa_{k}$ describes a pictorial element such as a train
or grass or a cow; $\randlabfa_{k} \in\indfa_{n}$ indicates that
picture $n$ contains, for example, a train; finally, the observed array
of pixels $X_{n}$
that form the picture is generated to contain the pictorial elements in
$\indfa_{n}$.
As in the clustering case, $X_{n}$ might not even be directly observed
but might
serve as a random effect in a deeper hierarchical model.

We observe that although the
order-of-appearance label sets are not exchangeable,
the sequence $(\indfa_{n})$ is. This fact
allows the formulation of a Gibbs sampler via the observation
that the distribution of $\indfa_{n}$, given the remaining
elements $\mathbf{\indfa}_{-n} := (\indfa_{1},\ldots,\indfa_{n-1},\break
\indfa_{n+1},\ldots, \indfa_{N})$, is the same as
if we thought of $\indfa_{n}$ as the final, $N$th element in a sequence
with $N-1$ preceding values given by $\mathbf{\indfa}_{-n}$.
The full details of such a sampler for the case of the IBP
(Examples~\ref{ex:epf_ibp}
and~\ref{ex:cond_ibp})
are given by \citet{griffiths:2006:infinite}.

As in the partition case, in practice, when programming the
sampler, it is useful to separate the feature allocation encoding from
the feature parameter values. \citet{griffiths:2006:infinite} describe how
\textit{left order form} matrices give a convenient representation
of the feature allocation in this context.

\section{Stick Lengths} \label{sec:stick}

Not every symmetric function defined for an arbitrary number
of arguments with values in the unit interval is an
EPPF \citep{pitman:1995:exchangeable}, and not every symmetric
function with an additional positive integer argument is an EFPF.
For instance, the consistency property in equation (\ref
{eq:strong_consistency})
implies certain additivity requirements for the function $\eppf$.

%
%
\begin{example}[(Not an EPPF)]
Consider the function~$p$ defined with
%
%
\begin{equation}
\label{eq:fake_eppf} \qquad p(1) = 1,\quad p(1,1) = 0.1,\quad p(2) = 0.8, \ldots
\end{equation}
From the information in equation (\ref{eq:fake_eppf}), $p$ may be
further defined so
as to be symmetric in its arguments
for any number of arguments,
but since it does not satisfy $p(1) = p(1,1) + p(2)$, it cannot be an EPPF.
\end{example}
%
%

%
%
\begin{example}[(Not an EFPF)]
Consider the function $p$ defined with
%
%
\begin{eqnarray}
\label{eq:fake_efpf} p(N=1) &= &0.9, \quad p(N=1,1) = 0.9,
\nonumber
\\[-8pt]
\\[-8pt]
\nonumber
\qquad p(N=1,1,1)
&=&
0.9, \ldots
\end{eqnarray}
From the information in equation (\ref{eq:fake_efpf}), $p$ may be
further defined so
as to be symmetric in its arguments
for any number of arguments after the initial $N$ argument,
but since $p(N=1) + p(N=1,1) + p(N=1,1,1) > 1$, it cannot be an
EFPF.
\end{example}
%
%

%
\begin{figure*}

\includegraphics{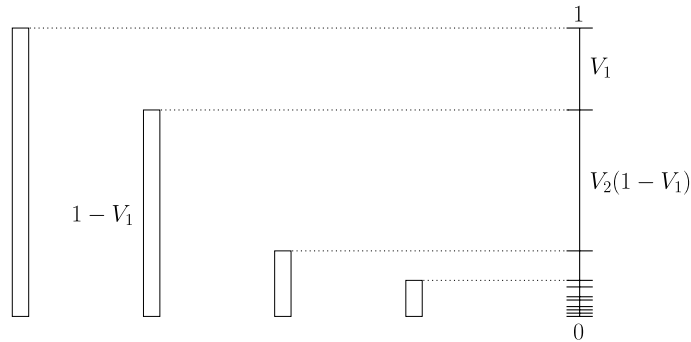}

\caption{An illustration of how stick-breaking
divides the unit interval into a sequence of probabilities Broderick, Jordan and
Pitman (\citeyear
{broderick:2012:beta}). The stick proportions $(\stickprop
_{1},\stickprop
_{2},\ldots)$ determine what fraction of the remaining stick is
appended to the probability sequence at each round.}\label{fig:stick_illus}
\end{figure*}

It therefore requires some care to define a suitable distribution
over consistent, exchangeable random feature allocations
or partitions using the exchangeable probability function
framework.

Since we are working with exchangeable sequences of random
variables, it is natural to turn to de Finet\-ti's
theorem (\cite{definetti1931funzione}; \cite{hewitt:1955:symmetric})
for clues as
to how to proceed. De Finetti's theorem tells us that any exchangeable
sequence of random variables can be expressed as an
independent and identically distributed sequence when conditioned
on an underlying random \textit{mixing measure}. While this theorem
may seem difficult to apply directly to, for example, exchangeable partitions,
it may be applied more naturally to an exchangeable
sequence of numbers derived from a sequence of partitions.
The argument below is
due to \citet{aldous:1985:exchangeability}.

Suppose that $(\randpart_{n})$ is an exchangeable, consistent
sequence of random partitions. Consider the $k$th partition block
to appear according to an order-of-appearance labeling scheme, and
give this block a new \textit{random} label, $\randlabpart_{k} \sim
\unif([0,1])$,
such that each random label is drawn independently from the rest.
This construction is the same as the one used for parameter generation
in Section~\ref{sec:epf_infer}, and $(\randpart_{n})$ is exchangeable
by the same
arguments used there. Let $\indpart_{n}$ equal $\randlabpart_{k}$
exactly when
$n$ belongs to the partition with this label.

If we apply de Finetti's theorem to the sequence $(\indpart_{n})$ and note
that $(\indpart_{n})$ has at most countably many different values, we
see that
there exists some random sequence $(\partblockfreq_{k})$ such that
$\partblockfreq_{k} \in(0,1]$ for all $k$ and, conditioned on the frequencies
$(\partblockfreq_{k})$, $(\indpart_{n})$ has the same distribution as
i.i.d. draws from
$(\partblockfreq_{k})$. In this description, we have brushed over
technicalities associated with partition blocks that contain only one
index even
as $N \rightarrow\infty$ (which may imply $\sum_{k} \partblockfreq_{k}
< 1$).

But
if we assume that every partition block eventually contains at least
two indices,
we can achieve an exchangeable partition of $[N]$ as follows.
Let $(\partblockfreq_{k})$ represent a sequence of values in $(0,1]$
such that
$\sum_{k=1}^{\infty} \partblockfreq_{k} \stackrel{\mathrm{a.s.}}{=} 1$. Draw
$\indpart_{n} \iid\disc((\partblockfreq_{k})_{k})$.
Let $\Pi_{N}$ be the induced partition given $(\indpart
_{n})_{n=1}^{N}$. Exchangeability
follows from the i.i.d. draws, and consistency follows from the induced
partition construction.

When the frequencies $(\partblockfreq_{k})$ are thought of as
subintervals of the unit interval,
that is, a partition of the unit interval, they are collectively
called \textit{Kingman's paintbox} \citep{kingman:1978:representation}.
As another naming convention, we may think of the unit interval as a
\textit{stick} \citep{ishwaran:2001:gibbs}. We
partition the unit interval by breaking it into various \textit{stick
lengths}, which
represent the frequencies of each partition block.

A similar construction can be seen to yield exchangeable, consistent
random feature allocations.
In this case, let $(\fablockfreq_{k})$ represent a sequence of values
in $(0,1]$ such that
$\sum_{k=1}^{\infty} \fablockfreq_{k} \stackrel{\mathrm{a.s.}}{<} \infty$. We
generate feature collections
independently for each index as follows. Start with $\indfa_{n} =
\varnothing$. For each feature $k$, add $k$ to the set $\indfa_{n}$,
independently from all other features, with probability $\fablockfreq_{k}$.
Let $\randfa_{N}$ be the induced feature allocation given $(\indfa
_{n})_{n=1}^{N}$.
Exchangeability of $\randfa_{N}$ follows from the i.i.d. draws of $\indfa
_{n}$, and consistency
follows from the induced feature allocation construction. The finite
sum constraint ensures each index
belongs to a finite number of features a.s.

It remains to specify a distribution on the partition or feature
frequencies. The frequencies cannot be
i.i.d. due to the finite summation constraint in both cases. In the
partition case, any infinite set of frequencies cannot even be
independent since the summation is fixed to one. One scheme to ensure
summation to unity is called
\textit{stick-breaking} (\cite{mccloskey:1965:model}; \cite{patil:1977:diversity};
\cite{sethuraman:1994:constructive}; \cite{ishwaran:2001:gibbs}). In stick-breaking,
the stick lengths are obtained by
recursively
breaking off parts of the unit interval to return as the
atoms $\partblockfreq_{1}, \partblockfreq_{2}, \ldots$ (cf.
Figure~\ref{fig:stick_illus}). In particular, we generate
stick-breaking proportions
$\stickprop_{1}, \stickprop_{2}, \ldots$ as $[0,1]$-valued random variables.
Then $\partblockfreq_{1}$ is the first proportion $\stickprop_{1}$
times the initial
stick length $1$; hence, $\partblockfreq_{1} = \stickprop_{1}$.
Recursively, after
$k$ breaks, the remaining length of the initial unit
interval is $\prod_{j=1}^{k} (1-\stickprop_{j})$. And $\partblockfreq
_{k+1}$ is the
proportion $\stickprop_{k+1}$ of the remaining stick; hence,
$\partblockfreq_{k+1} = \stickprop_{k+1} \prod_{j=1}^{k}
(1-\stickprop_{j})$.

The stick-breaking construction
yields $\partblockfreq_{1},\partblockfreq_{2},\ldots$ such that
$\partblockfreq_{k} \in[0,1]$ for each $k$ and
$\sum_{k=1}^{\infty} \partblockfreq_{k} \le1$. If the $\stickprop_{k}$
do not decay too rapidly, we will have $\sum_{k=1}^{\infty}
\partblockfreq_{k} \stackrel{\mathrm{a.s.}}{=} 1$.
In particular, the partition block proportions $\partblockfreq_{k}$
sum to unity a.s. iff there is no remaining stick
mass: $\prod_{k=1}^{\infty} (1-\stickprop_{k}) \stackrel{\mathrm{a.s.}}{=} 0$.

We often make the additional, convenient assumption that the
$\stickprop_{k}$ are independent. In this case, a
necessary and sufficient condition for $\sum_{k=1}^{\infty}
\partblockfreq_{k} \stackrel{\mathrm{a.s.}}{=} 1$
is $\sum_{k=1}^{\infty} \mbe [ \log(1-\stickprop_{k})
] =
-\infty$ \citep{ishwaran:2001:gibbs}.
When the $\stickprop_{k}$ are independent and of a canonical distribution,
they are easily simulated. Moreover,
if we assume that the $\stickprop_{k}$ are such that the
$\partblockfreq_{k}$
decay sufficiently rapidly in $k$, one strategy for
simulating a stick-breaking model is to ignore all $k > K$
for some fixed, finite $K$. This approximation is known as
truncation \citep{ishwaran:2001:gibbs}. It is fortuitously the case
that in some models of particular interest, such useful assumptions
fall out naturally from the model construction (e.g., Examples~\ref
{ex:stick_crp}
and~\ref{ex:stick_ibp}).

%
%
\begin{example}[(Chinese restaurant process)] \label{ex:stick_crp}
In the original exchangeability result due to de Finetti \citep
{definetti1931funzione},
the exchangeable random variables were zero/one-valued, and the
mixing measure was a distribution on a single frequency so that the
outcomes were conditionally Bernoulli. We will find a similar result in
obtaining the stick-breaking proportions associated with
the Chinese restaurant process.

We can construct a sequence of binary-valued random variables by
dividing the customers in the CRP who are sitting at the first table
from the rest; color the former collection of customers gray and the
latter collection of customers white. Then, we see that the first
customer must be colored gray. And thus we begin with a single gray
customer and no white customers. This binary valuation for the first
table in the CRP is illustrated by the first column in the matrix in
Figure~\ref{fig:polya_dp}.

%
\begin{figure*}

\includegraphics{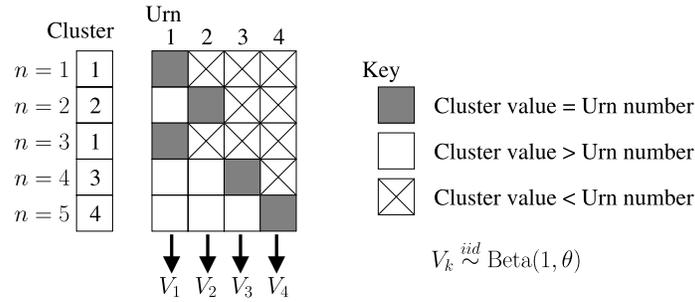}

\caption{An illustration of the proof based on
the P\'olya urn that Dirichlet process stick-breaking gives the
underlying partition block frequencies for a Chinese restaurant process
model. The $k$th column in the central matrix corresponds to a tallying
of when the $k$th table is chosen (gray), when a table of index larger
than $k$ is chosen (white), and when an index smaller than $k$ is
chosen ($\times$). If we ignore the $\times$ tallies, the gray and
white tallies in each column (after the first) can be modeled as balls
drawn from a P\'olya urn. The limiting frequency of gray balls in each
column is shown below the matrix.}\label{fig:polya_dp}
\end{figure*}

At this point, it is useful to recall the P\'olya urn construction
(\cite{polya:1930:sur};
\cite{freedman:1965:bernard}),\break whereby an urn starts with $G_{0}$ gray balls
and $W_{0}$ white balls. At each round $N$, we draw a ball from the
urn, replace it, and add $\urnextra$ of the same color of ball to the
urn. At the end of the round, we have $G_{N}$ gray balls and $W_{N}$
white balls. Despite the urn metaphor, the number of balls need not be
an integer at any time. By checking equation~(\ref{eq:crp}), which
defines the CRP, we
can see that the coloring of the gray/white customer matrix assignments
starting with the second customer has the same distributions as a
sequence of balls from a P\'olya urn as a P\'olya urn with $G_{1,0} =
1$ initial gray balls, $W_{1,0} = \dpconc$ initial white balls and
$\urnextra_{1} = 1$ replacement balls. Let $G_{1,N}$ and $W_{1,N}$
represent the numbers of gray and white balls, respectively, in the urn
after $N$ rounds. The important fact about the P\'olya urn we use here
is that there exists some $\stickprop\sim\tb(G_{0}/\urnextra,
W_{0}/\urnextra)$ such that $\urnextra^{-1}(G_{N+1}-G_{N}) \stackrel
{\mathrm{i.i.d.}}{\sim} \bern(\stickprop)$ for all $N$. In this particular case of
the CRP, then, $G_{1,N+1} - G_{1,N}$ is one if a customer sits at the
first table (or zero otherwise), and $G_{1,N+1} - G_{1,N} \stackrel
{\mathrm{i.i.d.}}{\sim} \bern(\stickprop_{1})$ with $\stickprop_{1} \sim\tb
(1,\dpconc)$.

We now look at the sequence of customers who sit at the second and
subsequent tables. That is, we condition on customers not sitting at
the first table or equivalently on the sequence with $G_{1,N+1} -
G_{1,N} = 0$. Again, we have that the first customer sits at the second
table, by the CRP construction. Now let customers at the second table
be colored gray and customers at the third and later tables be colored
white. This valuation is illustrated in the second column in
Figure~\ref{fig:polya_dp}; each $\times$ in the figure denotes a data
point where the
first partition block is chosen and, therefore, the current P\'olya urn
is not in play. As before, we begin with one gray customer and no white
customers. We can check equation (\ref{eq:crp}) to see that customer
coloring once
more proceeds according to a P\'olya urn scheme with $G_{2,0} = 1$
initial gray balls, $W_{2,0} = \dpconc$ initial white balls and
$\urnextra_{2} = 1$ replacement balls. Thus, contingent on a customer
not sitting at the first table, the $N$th customer sits at the second
table with i.i.d. distribution $\bern(\stickprop_{2})$ with $\stickprop
_{2} \sim\tb(1,\dpconc)$. Since the sequence of individuals sitting at
the second table has no other dependence on the sequence of individuals
sitting at the first table, we have that $\stickprop_{2}$ is
independent of $\stickprop_{1}$.

The argument just outlined proceeds recursively to show us that the
$N$th customer, conditional on not sitting at the first $K-1$ tables
for $K \ge1$, sits at the $K$th table with i.i.d. distribution $\bern
(\stickprop_{K})$ and $\stickprop_{K} \sim\tb(1,\dpconc)$ with
$\stickprop_{K}$ independent of the previous $(\stickprop_{1},\ldots
,\stickprop_{K-1})$.

%
\begin{figure*}

\includegraphics{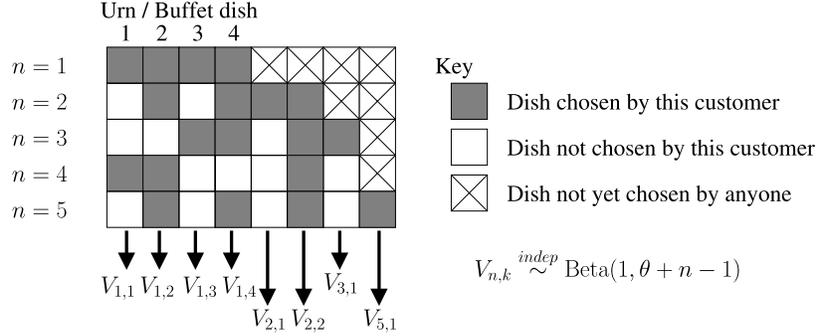}

\caption{Illustration of the proof that the
frequencies of features in the Indian buffet process are given by beta
random variables. For each feature, we can construct a sequence of
zero/one variables by tallying whether (gray, one) or not (white, zero)
that feature is represented by the given data point. Before the first
time a feature is chosen, we mark it with an $\times$. Each column
sequence of gray and white tallies, where we ignore the $\times$ marks,
forms a P\'olya urn with limiting frequencies shown below the
matrix.}\label{fig:polya_ibp}
\end{figure*}

Combining these results, we see that we have the following construction
for the customer seating patterns. The $\stickprop_{k}$ are distributed
independently and identically according to $\tb(1,\dpconc)$. The
probability $\partblockfreq_{K}$ of sitting at the $K$th table is the
probability of not sitting at the first $K-1$ tables, conditional on
not sitting at the previous table, times the conditional probability of
sitting at the $K$th table: $\partblockfreq_{K} =\break [ \prod_{k=1}^{K-1} (1-\stickprop_{k})  ] \cdot\stickprop_{K}$. Finally,
with the vector of table frequencies $(\partblockfreq_{k})$, each
customer sits independently and identically at the corresponding vector
of tables according to these frequencies. This process is summarized here:
%
%
\begin{eqnarray}\label{eq:dp_stick}
\nonumber
\stickprop_{k} &\stackrel{\mathrm{i.i.d.}} {\sim}& \tb(1,\dpconc),
\\
\partblockfreq_{K} &:=& \stickprop_{K} \prod
_{k=1}^{K} (1-\stickprop _{k}),
\\
 \indpart_{n} &\stackrel{\mathrm{i.i.d.}} {
\sim}& \disc\bigl((\partblockfreq_{k})_{k}\bigr).\nonumber
\end{eqnarray}

To see that this process is well-defined, first note that $\mbe [
\log(1-\stickprop_{k})  ]$ exists, is negative and is the same for
all $k$ values. It follows that $\sum_{k=1}^{\infty} \mbe [
\log
(1-\stickprop_{k})  ] = -\infty$, so by the discussion before this
example, we must have $\sum_{k=1}^{K} \partblockfreq_{k} \stackrel
{\mathrm{a.s.}}{=} 1$.
\end{example}

The feature case is easier. Since it does not require the frequencies
to sum to one, the random frequencies can be independent so long as
they have an a.s. finite sum.

%
%
\begin{example}[(Indian buffet process)] \label{ex:stick_ibp}
As in the case of the CRP, we can recover the stick lengths for the
Indian buffet
process using an argument based on an urn model.

Recall that on the first round of the Indian buffet process, $K_{1}^{+}
\sim\pois(\bpmass)$ features are chosen to contain index $1$. Consider
one of the features, labeled~$k$. By construction, each future data
point $N$ belongs to this feature with probability $N_{N-1,k} /
(\bpconc+ N -\nobreak 1)$. Thus, we can model the sequence after the first
data point as a P\'olya urn of the sort encountered in Example \ref
{ex:stick_crp}
with initially $G_{k,0} = 1$ gray balls, $W_{k,0} = \bpconc$ white
balls and $\urnextra_{k} = 1$ replacement balls. As we have seen, there
exists a random variable $V_{k} \sim\tb(1,\bpconc)$ such that
representation of this feature by data point $N$ is chosen, i.i.d. across
all $N$, as $\bern(\stickprop_{k})$. Since the Bernoulli draws
conditional on previous draws are independent across all $k$, the
$\stickprop_{k}$ are likewise independent of each other; this fact is
also true for $k$ in future rounds.\vadjust{\goodbreak} Draws according to such an urn are
illustrated in each of the first four columns of the matrix in
Figure~\ref{fig:polya_ibp}.

Now consider any round $n$. According to the IBP construction,
$K_{n}^{+} \sim\pois(\bpmass\bpconc/ (\bpconc+ n - 1))$ new
features are chosen to include index $n$. Each future data point $N$
(with $N > n$) represents feature $k$ among these features with
probability $N_{N-1,k} / (\bpconc+ N - 1)$. In this case, we can model
the sequence after the $n$th data point as a P\'olya urn with $G_{k,0}
= 1$ initial gray balls, $W_{k,0} = \bpconc+ n - 1$ initial white
balls and $\urnextra_{k} = 1$ replacement balls. So there exists a
random variable $\stickprop_{k} \sim\tb(1,\bpconc+ n - 1)$ such that
representation of feature $k$ by data point $N$ is chosen, i.i.d. across
all $N$, as $\bern(\stickprop_{k})$.

Finally, then, we have the following generative model for the feature
allocation by iterating across $n=1,\ldots,N$ \citep{ThibauxJo07}:
%
%
\begin{eqnarray}
\label{eq:beta_proc_num_sticks}
K_{n}^{+} &\stackrel{\mathrm{indep}} {\sim}& \pois \biggl(
\frac{\bpmass
\bpconc
}{\bpconc+ n - 1} \biggr),
\\
\nonumber
K_{n} &= &K_{n-1} + K_{n}^{+},
\\
\label{eq:beta_proc_stick_lengths}
\stickprop_{k} &\stackrel{\mathrm{indep}} {\sim}& \tb(1,\bpconc+ n - 1),
\nonumber
\\[-8pt]
\\[-8pt]
\eqntext{ k =
K_{n-1} + 1,\ldots,K_{n},}
\\
\nonumber
\indhasf_{n,k} &\stackrel{\mathrm{indep}} {\sim}& \bern(
\stickprop_{k}),\quad  k = 1, \ldots, K_{n}.
\end{eqnarray}
$\indhasf_{n,k}$ is an indicator random variable for whether feature
$k$ contains index $n$. The collection of features to which index $n$
belongs, $\indfa_{n}$, is the collection of features $k$ with
$\indhasf
_{n,k} = 1$.
\end{example}
%
%

\subsection{Inference} \label{sec:stick_infer}

As we have seen above, the exchangeable probability functions of
Section~\ref{sec:epf}
are the marginal distributions of the partitions or feature allocations
generated according to stick-length models with the stick lengths
integrated out. It has been proposed that including the stick lengths
in MCMC samplers of these models will improve mixing (Ishwaran and
Zarepour, \citeyear{ishwaran:2000:markov}). While it
is impossible to sample the countably infinite set of partition block
or feature frequencies in these models (cf. Examples \ref
{ex:stick_crp} and~\ref{ex:stick_ibp}), a~number
of ways of getting around this difficulty have been investigated.
\citet{ishwaran:2000:markov} examine two separate finite approximations
to the full CRP
stick-length model: one uses a parametric approximation to the full
infinite model, and the other creates a truncation by setting the stick
break at some fixed size $K$ to be~1: $\stickprop_{K} = 1$.
There also exist techniques that avoid any approximations and deal
instead directly with the full model, in particular,
retrospective sampling (Papaspiliopoulos and Roberts, \citeyear{papaspiliopoulos:2008:retrospective}) and
slice sampling \citep{walker:2007:sampling}.

While our discussion thus far has focused on\break MCMC
sampling as a means of approximating the posterior distribution
of either the block assignments or both the block assignments and
stick lengths, including the stick lengths in a posterior analysis
facilitates a different posterior
approximation; in particular, \textit{variational methods} can also be
used to approximate the posterior. These methods minimize some notion of
distance to the posterior over a family of
potential approximating distributions \citep{jordan:1999:introduction}.
The practicality and, indeed,
speed of these methods in the case of stick-breaking for the CRP
(Example~\ref{ex:stick_crp}) have been
demonstrated by \citet{blei:2006:variational}.

A number of different models for the stick lengths corresponding
to the features of an IBP (Example~\ref{ex:stick_ibp}) have been
discovered. The
distributions described in Example~\ref{ex:stick_ibp} are covered by
\citet{ThibauxJo07}, who build on work from
\citet{hjort:1990:nonparametric}, \citet{Kim99}. A
special case of the IBP is examined by \citet{teh:2007:stick}, who
detail a slice sampling algorithm for sampling from the
posterior of the stick lengths and feature assignments.
Yet another stick-length model for the IBP is explored by
\citet{paisley:2010:stick}, who show how to apply variational methods to
approximate the posterior of their model.

Stick-length modeling has the further advantage of allowing
inference in cases where it is not straightforward to integrate
out the underlying stick lengths to obtain a tractable exchangeable
probability function.

\section{Subordinators} \label{sec:sub}

An important point to reiterate about the labels $\indpart_{n}$ and
label collections $\indfa_{n}$
is that when we use the order-of-appearance labeling scheme for
partition or feature blocks
described above, the random sequences $(\indpart_{n})$ and $(\indfa
_{n})$ are not exchangeable.
Often, however, we would like to make use of special properties of
exchangeability
when dealing with these sequences. For instance, if we use Markov Chain
Monte Carlo
to sample from the posterior distribution of a partition (cf.
Section~\ref{sec:epf_infer}), we might want to Gibbs sample the
cluster assignment of
data point $n$ given the assignments of the remaining data points:
$\indpart_{n}$ given $\{\indpart_{m}\}_{m=1}^{N} \setminus\{
\indpart
_{n}\}$. This sampling is particularly
easy in some cases \citep{neal:2000:markov} if we can treat $\indpart
_{n}$ as the last random
variable in the sequence, but this treatment requires exchangeability.

A way to get around this dilemma was suggested by \citet
{aldous:1985:exchangeability} and appeared above in
our motivation for using stick lengths. Namely, we assign to the $k$th
partition block a
uniform random label $\randlabpart_{k} \sim\unif([0,1])$; analogously,
we assign to the $k$th
feature a uniform random label $\randlabfa_{k} \sim\unif([0,1])$.
We can see that in both cases, all of the labels are a.s. distinct.
Now, in the partition case, let $\indpart_{n}$
be the uniform random label of the partition block to which $n$
belongs. And in the feature case,
let $\indfa_{n}$ be the (finite) set of uniform random feature labels
for the
features to which $n$ belongs. We can recover the partition or feature
allocation as the
induced partition or feature allocation by grouping indices assigned to
the same label. Moreover, as discussed above, we now have that each of
$(\indpart_{n})$ and
$(\indfa_{n})$ is an exchangeable sequence.

%
\begin{figure*}

\includegraphics{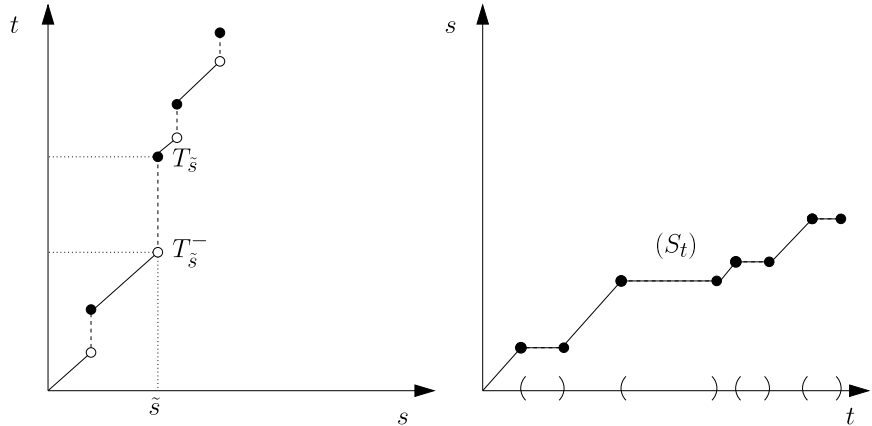}

\caption{\textit{Left}: The sample path $(\subord
_{s})$ of a subordinator. $\subord_{\tilde{s}}^{-}$ is the limit from
the left of $(\subord_{s})$ at $s = \tilde{s}$. \textit{Right}: The
right-continuous inverse $(\invsub_{t})$ of a subordinator: $S_{t} :=
\inf\{s\dvtx T_{s} > t\}$. The open intervals along the $t$ axis correspond
to the jumps of the subordinator $(\subord_{s})$.}\label{fig:subordinator}
\end{figure*}

If we form partitions or features according to the stick-length
constructions detailed in Section~\ref{sec:stick}, we know that each unique
partition or feature label $\rmloc_{k}$ is associated with a frequency
$\rmweight_{k}$.
We can use this association to form a random measure:
%
%
\begin{equation}
\label{eq:rand_meas} \mu= \sum_{k=1}^{\infty}
\rmweight_{k} \delta_{\rmloc_{k}},
\end{equation}
where $\delta_{\rmloc_{k}}$ is a unit point mass located at $\rmloc_{k}$.
In the partition case, $\sum_{k} \rmweight_{k} = 1$,\vspace*{-1pt} so the random measure
is a random probability measure, and we may draw $\indpart_{n} \iid
\mu$.
In the feature case, the weights have a finite sum
but do not necessarily sum to one. In the feature case, we draw $\indfa
_{n}$ by including each
$\rmloc_{k}$ for which $\bern(\rmweight_{k})$ yields a draw of 1.

Another way to codify the random measure in equation (\ref{eq:rand_meas})
is as a monotone increasing stochastic process on $[0,1]$. Let
\[
\subord_{s} = \sum_{k=1}^{\infty}
\rmweight_{k} \mbo\{\rmloc_{k} \le s\}.
\]
Then the atoms of $\mu$ are in one-to-one correspondence with the jumps
of the process $\subord$.

This increasing random function construction gives us another means of
choosing distributions for the weights $\rmweight_{k}$. We have already
seen that
these cannot be i.i.d. due to the finite summation condition. However, we will
see that if we require that the \textit{increments} of a monotone, increasing
stochastic process are independent and stationary, then we can use the
jumps of that function as the atoms in our random measure for
partitions or features.

\begin{definition} $\!\!$A \textit{subordinator}
(\cite{bochner:1955:harmonic};
\citeauthor{bertoin:1998:levy} \citeyear{bertoin:1998:levy,bertoin2004:subordinatorsEX})
is a stochastic process $(\subord_{s}, s \ge0)$ that has the
following properties:
\begin{itemize}
\item Nonnegative, nondecreasing paths (a.s.),
\item Paths that are right-continuous with left limits, and
\item Stationary, independent increments.
\end{itemize}
\end{definition}
For our purposes, wherein the subordinator values will ultimately
correspond to (perhaps scaled) probabilities, we will assume the
subordinator takes values in $[0, \infty)$, though alternative ranges
with a sense of ordering are possible.

Subordinators are of interest to us because they not only exhibit
the stationary independent increments property\vadjust{\goodbreak} but they also can
always be decomposed into two components: a deterministic \textit{drift}
component and a \textit{Poisson point process}. Recall that a Poisson
point process on space $S$ with rate measure $\intens(d x)$, where $x
\in S$, yields a countable subset of points of $S$. Let $N(A)$ be the
number of points of the process in set $A$ for $A \subseteq S$. The
process is characterized by the fact that, first, $N(A) \sim\pois
(\intens(A))$ for any $A$ and, second, for any disjoint $A_{1}, \ldots,
A_{K}$, we have that $N(A_{1}), \ldots, N(A_{K})$ are independent
random variables. See \citet{kingman:1993:poisson} for a thorough
treatment of these processes. An example subordinator with both drift
and jump components is shown on the left-hand side of Figure~\ref{fig:subordinator}.

The subordinator decomposition is detailed in the following result
\citep{bertoin:1998:levy}.
%
\begin{theorem} \label{thm:sub_ppp}
Every subordinator $(\subord_{s}, s \ge0)$ can be written as
%
%
\begin{equation}
\label{eq:sub_drift_jump} \subord_{s} =
\drift s + \sum_{k=1}^{\infty}
\subjumpsize_{k} \mbo \{ \subjumploc_{k} \le s\}
\end{equation}
for some constant $\drift\ge0$ and where $\{(\subjumpsize
_{k},\subjumploc_{k})\}_{k}$ is the countable set of points of a
Poisson point process with intensity $\levym(d\subjumpsize)
\,d\subjumploc$, where $\levym$ is a L\'{e}vy measure; that is,
\[
\int_{0}^{\infty} (1 \wedge\subjumpsize) \levym(d
\subjumpsize) < \infty.
\]
\end{theorem}
%
In particular, then, if a subordinator is finite at time~$t$, the jumps
of the subordinator up to $t$ may be used as
feature block frequencies if they have support in $[0,1]$. Or, in
general, the normalized jumps may be used as partition block frequencies.
We can see from the right-hand side of Figure~\ref{fig:subordinator}
that the
jumps of a subordinator partition intervals of the form $[0,t)$, as
long as the subordinator has no drift component.
In either the feature or cluster case,
we have substituted the condition of independent and identical
distribution for the partition or feature frequencies (i.e., the jumps)
with a more natural continuous-time analogue: independent, stationary intervals.

Just as the Laplace transform of a positive random variable
characterizes the distribution of that random variable,
so does the Laplace transform of the subordinator---which is a positive
random variable at
any fixed time point---describe this stochastic process (\citeauthor{bertoin:1998:levy}
\citeyear
{bertoin:1998:levy,bertoin2004:subordinatorsEX}).
%
\begin{theorem}[(L\'{e}vy--Khinchin formula for subordinators)]
\label{thm:lk}
If $(\subord_{s}, s \ge0)$ is a subordinator, then for $\lambda\ge0$
we have
%
%
\begin{equation}
\label{eq:laplace_exp} \mathbb{E}\bigl(e^{-\lambda\subord_{s}}\bigr) =
e^{-\lapexp(\lambda) s}
\end{equation}
with
%
%
\begin{equation}
\label{eq:lk} \lapexp(\lambda) = \drift\lambda+ \int_{0}^{\infty}
\bigl(1-e^{-\lambda\subjumpsize} \bigr) \levym(d \subjumpsize),
\end{equation}
where $\drift\ge0$ is called the drift constant and $\levym$ is a
nonnegative, L\'{e}vy measure on $(0,\infty)$.
\end{theorem}
%
The function $\lapexp(\lambda)$ is called the \textit{Laplace
exponent} in
this context.
We note that a subordinator is characterized by its drift constant and
L\'{e}vy measure.

Using subordinators for feature allocation modeling is particularly
easy; since the
jumps of the subordinators are formed by a Poisson point process, we
can use Poisson
process methodology to find the stick lengths and EFPF. To set up this
derivation,
suppose we generate feature membership from a subordinator by taking Bernoulli
draws at each of its jumps with success probability equal to the jump size.
Since every jump has strictly positive size, the feature associated
with each
jump will eventually score a Bernoulli success for some index $n$ with
probability
one. Therefore, we can enumerate all jumps of the process in order of
appearance;
that is, we first enumerate
all features in which index $1$ appears, then all features in which
index $2$
appears but not index $1$, and so on. At the $n$th iteration, we
enumerate all
features in which index $n$ appears but not previous indices. Let $K_{n}^{+}$
represent the number of indices so chosen on the $n$th round. Let
$K_{0} = 0$
so that recursively $K_{n} := K_{n-1} + K_{n}^{+}$ is the number of subordinator
jumps seen by round $n$, inclusive. Let $\subjumpsize_{k}$ for
$k = K_{n-1} + 1,\ldots,K_{n}$ be the distribution of a particular subordinator
jump seen on round $n$. We now turn to connecting the subordinator perspective
to the earlier derivation of stick lengths in Section~\ref{sec:stick}.

%
%
\begin{example}[(Indian buffet process)] \label{ex:sub_ibp}
In our earlier discussion, we found a collection of stick lengths to
represent the featural frequencies for the IBP [equation~(\ref
{eq:beta_proc_stick_lengths}) of Example~\ref{ex:stick_ibp} in
Section~\ref{sec:stick}]. To see
the connection to subordinators, we start from the
\textit{beta process subordinator} \citep{Kim99} with zero drift
($\drift
= 0$) and L\'{e}vy measure
%
%
\begin{equation}
\label{eq:beta_levy} \levym(d \rmweight) = \bpmass\bpconc
\rmweight^{-1} (1-\rmweight )^{\bpconc- 1} \,d\rmweight.
\end{equation}
We will see that the mass parameter $\bpmass> 0$ and concentration
parameter $\bpconc> 0$
are the same as those introduced in Example~\ref{ex:epf_ibp} and
continued in Example~\ref{ex:stick_ibp}.

\begin{theorem}
\label{thm:sub_sticks}
Generate a feature allocation from a beta process subordinator with L\'
{e}vy measure given by equation (\ref{eq:beta_levy}).
Then the sequence of subordinator jumps $(\subjumpsize_{k})$, indexed
in order of appearance, has the same distribution as the sequence of
IBP stick lengths $(\stickprop_{k})$ described by equations (\ref
{eq:beta_proc_num_sticks}) and (\ref{eq:beta_proc_stick_lengths}).
\end{theorem}

%
\begin{figure}[b]

\includegraphics{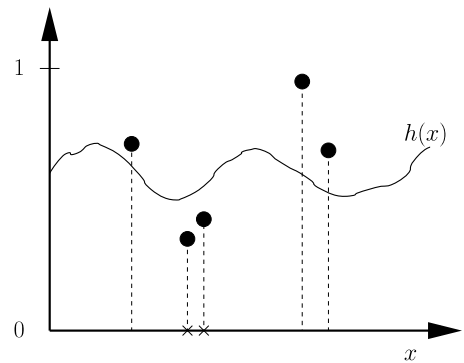}

\caption{An illustration of Poisson thinning. The
$x$-axis values of the filled black circles, emphasized by dotted
lines, are generated according to a Poisson process. The $[0,1]$-valued
function $h(x)$ is arbitrary. The vertical axis values of the points
are uniform draws in $[0,1]$. The ``thinned'' points are the collection
of $x$-axis values corresponding to vertical axis values below $h(x)$
and are denoted with a $\times$ symbol.}\label{fig:thinning}
\end{figure}

\begin{pf}
Recall the following fact about Poisson thinning \citep
{kingman:1993:poisson}, illustrated in Figure~\ref{fig:thinning}.
Suppose that a
Poisson point process with rate measure $\lambda$ generates points with
values $x$. Then suppose that, for each such point $x$, we keep it with
probability $h(x) \in[0,1]$. The resulting set of points is also a
Poisson point process, now with rate measure $\lambda'(A) = \int_{A}
\lambda(dx) h(x) \,dx$.


We prove Theorem~\ref{thm:sub_sticks} recursively. Define the measure
\[
\mu_{n}(d \rmweight) := \bpmass\bpconc\rmweight^{-1} (1-
\rmweight )^{\bpconc+ n - 1} \,d\rmweight,
\]
so that $\mu_{0}$ is the beta process L\'{e}vy measure $\levym$ in
equation (\ref{eq:beta_levy}).
We make the recursive assumption that $\mu_{n}$ is distributed as the
beta process
measure without atoms corresponding to features chosen on the first
$n$ iterations.

There are two parts to proving Theorem~\ref{thm:sub_sticks}. First, we
show that,
on the $n$th iteration,
the number of features chosen and the distribution of the corresponding
atom weights
agree with equations (\ref{eq:beta_proc_num_sticks}) and (\ref
{eq:beta_proc_stick_lengths}), respectively. Second, we check that the
recursion assumption
holds.

For the first part, note that on the $n$th round we choose features
with probability equal to their
atom weight. So we form a thinned Poisson process with rate measure
$\rmweight\cdot\mu_{n-1}(d \rmweight)$. This rate measure has total mass
\[
\int_{0}^{1} \rmweight\cdot\mu_{n-1}(d
\rmweight) = \bpmass\frac
{\bpconc}{\bpconc+ n - 1} =: \bpmass_{n-1}.
\]
So the number of features chosen is Poisson-distrib\-uted with mean
$\bpmass\bpconc(\bpconc+ n -1)^{-1}$, as desired [cf. equation~(\ref
{eq:beta_proc_num_sticks})]. And the atom weights have distribution equal
to the normalized rate measure
\begin{eqnarray*}
&&\bpmass_{n-1}^{-1} \rmweight\cdot\bpmass\bpconc
\rmweight^{-1} (1-\rmweight)^{\bpconc+ (n - 1) - 1} \,d\rmweight\\
&&\quad= \tb(\rmweight| 1,
\bpconc+ n - 1) \,d\rmweight
\end{eqnarray*}
as desired [cf. equation (\ref{eq:beta_proc_stick_lengths})].

Finally, to check the recursion assumption, we note that those sticks
that remain were chosen for having Bernoulli failure draws; that is,
they were chosen with probability equal to one minus their atom weight.
So the thinned rate measure for the next round is
\[
(1- \rmweight) \cdot\bpmass\bpconc\rmweight^{-1} (1-\rmweight
)^{\bpconc+ (n-1) - 1} \,d\rmweight,
\]
which is just $\mu_{n}$.
\end{pf}

The form of the EFPF of the feature allocation generated from the beta
process subordinator
follows immediately from the stick-length distributions we have just
derived by the discussion in Example~\ref{ex:stick_ibp} in
Section~\ref{sec:stick}.
\end{example}

We see from the previous example that feature allocation stick lengths
and EFPFs can be obtained in a straightforward manner\vadjust{\goodbreak} using the Poisson
process representation of the jumps of the subordinator. Partitions,
however, are not as easy to analyze, principally due to the fact that
the subordinator jumps must first be normalized to obtain a probability
measure on $[0,1]$; a random measure with finite total mass is not
sufficient in the partition case. Hence, we must compute the stick
lengths and EPPF using partition block frequencies from these
normalized jumps instead of directly from the subordinator jumps.

In the EPPF case, we make use of a result that gives us the
exchangeable probability function as a function of the Laplace
exponent. Though we do not derive this formula here, its derivation can
be found in \citet{pitman:2003:poisson}; the proof relies on, first,
calculating the joint distribution of the subordinator jumps and
partition generated from the normalized jumps and, second, integrating
out the subordinator jumps to find the partition marginal.
%
\begin{theorem}
\label{thm:sub_eppf}
Form a probability measure $\mu$ by normalizing jumps of the
subordinator with Laplace exponent $\lapexp$.
Let $(\randpart_{n})$ be a consistent set of exchangeable partitions
induced by i.i.d. draws from $\mu$.
For each exchangeable partition $\detpart_{N} = \{\blockpart
_{1},\ldots
,\blockpart_{K}\}$ of $[N]$ with
$N_{k} := |A_{k}|$ for each $k$,
%
%
\begin{eqnarray}\label{eq:laplace_eppf}
\quad&&\mathbb{P}(\randpart_{N} = \detpart_{N})\nonumber \\
&&\quad=
\eppf(N_{1},\ldots,N_{K})
\nonumber\\
&&\quad= \frac{(-1)^{N-K}}{(N-1)!} \int
_{0}^{\infty} \lambda^{N-1} e^{-\lapexp(\lambda)}
\prod_{k=1}^{K} \lapexp^{(N_{k})}(
\lambda) \,d\lambda,
\end{eqnarray}
where $\lapexp^{(N_{k})}(\lambda)$ is the $N_{k}$th derivative of the
Laplace exponent $\lapexp$ evaluated at $\lambda$.
\end{theorem}
%

%
%
\begin{example}[(Chinese restaurant process)] \label{ex:sub_crp}
We start by introducing the \textit{gamma process}, a subordinator
that we will see below generates the Chinese restaurant process EPPF.
The gamma process has Laplace exponent $\lapexp(\lambda)$
[equation (\ref{eq:laplace_exp})] characterized by
%
%
\begin{equation}
\label{eq:gamma_levy} \drift= 0\quad \mathrm{and}\quad \levym(d \rmweight) =
\gpconc \rmweight^{-1} e^{-\gpscale\rmweight} \,d \rmweight
\end{equation}
for $\gpconc> 0$ and $\gpscale> 0$ [cf. equation (\ref{eq:lk}) in
Theorem~\ref{thm:lk}]. We
will see that $\gpconc$ corresponds to the CRP concentration parameter
and that $\gpscale$ is arbitrary and does not affect the partition model.

We calculate the EPPF using Theorem~\ref{thm:sub_eppf}.
%
\begin{theorem} \label{thm:eppf_gamma}$\!\!\!\!$
The EPPF for partition block membership chosen according to the
normalized jumps $(\partblockfreq_{k})$ of the gamma subordinator with
parameter $\theta$ is the CRP EPPF [equation (\ref{eq:crp_eppf})].
\end{theorem}
%

\begin{pf}
By Theorem~\ref{thm:sub_eppf}, if we can find all order derivatives of
the Laplace
exponent $\lapexp$, we can calculate the EPPF for the partitions
generated with frequencies equal to the normalized jumps of this
subordinator. The derivatives of $\lapexp$, which are known to always
exist (\cite{bertoin2000:subordinatorsLE};
\cite{rogers:2000:diffusions}), are
straightforward to calculate if we begin by noting that, from
equation (\ref{eq:lk})
in Theorem~\ref{thm:lk}, we have in general that
\[
\lapexp'(\lambda) = \drift+ \int_{0}^{\infty}
\rmweight e^{-\lambda
\rmweight} \Lambda(d \rmweight).
\]
Hence, for the gamma process subordinator,
\[
\lapexp'(\lambda) = \int_{0}^{\infty}
e^{-\lambda\rmweight} \gpconc e^{-\gpscale\rmweight} \,d \rmweight= \frac{\gpconc}{\lambda+
\gpscale}.
\]
Then simple integration and differentiation yield
\[
\lapexp(\lambda) = \gpconc\log(\lambda+ \gpscale) - \gpconc\log (\gpscale)
\]
since $\lapexp(0) = 0$ and
\[
\lapexp^{(n)}(\lambda) = (-1)^{n-1} \frac{(n-1)! \gpconc}{(\lambda+
\gpscale)^{n}},\quad n
\ge1.
\]
We can substitute these quantities into the general EPPF formula in
equation (\ref{eq:laplace_eppf}) of Theorem~\ref{thm:sub_eppf} to
obtain
\begin{eqnarray*}
&& \eppf(N_{1},\ldots,N_{K})
\\
&&\quad= \frac{(-1)^{N-K}}{(N-1)!} \int_{0}^{\infty}
\lambda^{N-1} (\lambda + \gpscale)^{-\gpconc} \gpscale^{\gpconc}
\\
&&\hspace*{64pt}\qquad{}\cdot\prod_{k=1}^{K} (-1)^{N_{k}-1}
\frac{(N_{k} - 1)! \theta}{(\lambda+ \gpscale)^{N_{k}}} \,d\lambda
\\
&&\quad= \gpscale^{\gpconc} \frac{\gpconc^{K}}{(N-1)!} \Biggl[ \prod
_{k=1}^{K} (N_{k} - 1)! \Biggr]
\gpscale^{N-1-N-\gpconc+1} \\
&&\qquad{}\cdot\int_{0}^{\infty}
x^{N-1} (x + 1)^{-N-\gpconc} \,dx\quad\mathrm{for}\ x = \lambda/ \gpscale
\\
&&\quad= \frac{\gpconc^{K}}{(N-1)!} \Biggl[ \prod_{k=1}^{K}
(N_{k} - 1)! \Biggr] \frac{\Gamma(N) \Gamma(\gpconc)}{\Gamma(N+\gpconc)}
\\
&&\quad= \theta^{K} \Biggl[ \prod_{k=1}^{K}
(N_{k} - 1)! \Biggr] \frac
{1}{\gpconc(\gpconc+ 1)_{N-1 \uparrow1}}.
\end{eqnarray*}
The penultimate line follows from the form of the beta prime
distribution. The final line is the CRP EPPF from equation (\ref
{eq:crp_eppf}), as
desired. We note in particular that the parameter $\gpscale$ does not
appear in the final EPPF.
\end{pf}
\end{example}

Whenever the Laplace exponent of a subordinator is known,
Theorem~\ref{thm:sub_eppf} can similarly be applied to quickly find
the EPPF of the
partition generated by sampling from the normalized subordinator jumps.

To find the distributions of the stick lengths---that is, the partition
block frequencies---from the subordinator representation for a
partition, we must find the distributions of the normalized
subordinator jumps.

As in the feature case, we may enumerate the jumps of a subordinator
used for partitioning in the order of their appearance. That is, let
$\partblockfreq_{1}$
be the normalized subordinator jump size corresponding to the cluster
of the
first data point. Recursively, suppose index $n$ joins a cluster to
which none of the indices in $[n-1]$ belong, and suppose there are $k$ clusters
among $[n-1]$. Then let $\partblockfreq_{k+1}$ be the normalized subordinator
jump size corresponding to the cluster containing $n$.

%
%
\begin{example}[(Chinese restaurant process)] \label{ex:sub_crp_sticks}
We continue with the CRP example.
%
\begin{theorem} \label{thm:stick_gamma}$\!\!\!$
The normalized subordinator jumps $(\partblockfreq_{k})$ in order of
appearance of the gamma subordinator with concentration parameter
$\gpconc$ (and arbitrary parameter $\gpscale> 0$) have the same
distribution as the CRP stick lengths [equation (\ref{eq:dp_stick}) of
Example~\ref{ex:stick_crp}
in Section~\ref{sec:stick}].
\end{theorem}
%

\begin{pf}
First, we introduce some notation. Let $\sumofjumps= \sum_{k}
\subjumpsize_{k}$, the sum over all of the jumps of the subordinator.
Second, let
$\sumofjumps_{k} = \sumofjumps- \sum_{j=1}^{k} \subjumpsize_{k}$, the
total sum minus the first $k$ elements (in order of appearance). Note
that $\tau= \tau_{0}$. Finally, let $\omsprop_{k} = \sumofjumps_{k} /
\sumofjumps_{k-1}$ and $\stickprop_{k} = 1 - \omsprop_{k}$. Then a
simple telescoping of factors shows that $\partblockfreq_{k} =
\stickprop_{k} \prod_{j=1}^{k-1} (1-\stickprop_{j})$:
\begin{eqnarray*}
\stickprop_{k} \prod_{j=1}^{k-1}
(1-\stickprop_{j})& = &\biggl(1 - \frac{\sumofjumps_{k}}{\sumofjumps_{k-1}} \biggr) \prod
_{j=1}^{k-1} \frac{\sumofjumps_{j}}{\sumofjumps_{j-1}} \\
&=&
\frac{\sumofjumps_{k-1} - \sumofjumps_{k}}{\sumofjumps_{0}} = \frac{\subjumpsize_{k}}{\sumofjumps} = \partblockfreq_{k}.
\end{eqnarray*}

It remains to show that the $\stickprop_{k}$ have the desired
distribution. To that end, it is easier to work with the $\omsprop
_{k}$. We will find the following lemma \citep
{pitman:2006:combinatorial} useful.
%
\begin{lemma}
Consider a subordinator with L\'{e}vy measure $\levym$, and suppose
$\sumofjumps$ equals the sum of all jumps of the subordinator.
Let $\levydens$ be the density of $\levym$ with respect to the Lebesgue
measure. And let $\sumdens$ be the density of\vadjust{\goodbreak} the distribution of
$\sumofjumps$ with respect to the Lebesgue measure. Then
\begin{eqnarray*}
&& \mbp(\sumofjumps_{0} \in dt_{0}, \ldots,
\sumofjumps_{k} \in dt_{k})
\\
&&\quad= f(t_{k}) \,dt_{k} \Biggl( \prod
_{j=0}^{k-1} \frac{(t_{j} -
t_{j+1})\rho(t_{j} - t_{j+1})}{t_{j}} \,dt_{j}
\Biggr).
\end{eqnarray*}
\end{lemma}
%

With this lemma in hand, the result follows from a change of variables
calculation; we use a bijection between $\{\omsprop_{1},\ldots
,\omsprop
_{k},\tau\}$ and $\{\sumofjumps_{0},\ldots,\sumofjumps_{k}\}$ defined
by $\sumofjumps_{k} = \sumofjumps\prod_{j=1}^{k} \omsprop_{j}$. The
determinant of the Jacobian for the transformation to the former
variables from the latter is
\[
J = \prod_{j=1}^{k} \Biggl[ \tau\prod
_{i=1}^{j-1} W_{i} \Biggr] =
\prod_{j=0}^{k-1} \sumofjumps_{j}(
\tau, W_{1}, \ldots, W_{k}).
\]
In the derivation that follows, we start by expressing results in terms
of the $\sumofjumps_{j}$ terms with the dependence on $\{\tau, W_{1},
\ldots, W_{k}\}$ suppressed to avoid notational clutter, for example,
$J = \prod_{j=0}^{k-1} \sumofjumps_{j}$. At the end, we will evaluate
the $\sumofjumps_{j}$ terms as functions of $\{\tau, W_{1}, \ldots,
W_{k}\}$.

For now, then, we have
\begin{eqnarray*}
&& \mathbb{P}(\omsprop_{1} \in d\soms_{1}, \ldots,
\omsprop_{k} \in d\soms_{k}, \tau\in dt_{0})
\\
&&\quad= \mathbb{P}(\sumofjumps_{0} \in dt_{0}, \ldots,
\sumofjumps_{k} \in dt_{k}) \cdot J
\\
&&\quad= \sumdens(t_{k}) \,dt_{k} \Biggl( \prod
_{j=0}^{k-1} (t_{j} - t_{j+1})
\levydens(t_{j} - t_{j+1}) \Biggr).
\end{eqnarray*}

In the case of the gamma process, we can read $\levydens( \subjumpsize
) = \gpconc\subjumpsize^{-1} e^{-\gpscale\subjumpsize}$ from
equation (\ref{eq:gamma_levy}). The function $\sumdens$ is determined
by $\levydens$ and
in this case \citep{pitman:2006:combinatorial},
\[
\sumdens(t) = \ga(t | \gpconc, \gpscale) = \gpscale^{\gpconc} \Gamma (
\gpconc)^{-1} t^{\gpconc- 1} e^{-\gpscale t}.
\]
So
\begin{eqnarray*}
&& \mathbb{P}(\omsprop_{1} \in d\soms_{1}, \ldots,
\omsprop_{k} \in d\soms_{k}, \sumofjumps\in
dt_{0})\\
&&\quad\propto t_{k}^{\gpconc- 1} e^{-\gpscale t_{0}} = t_{0}^{\gpconc-1}
e^{-\gpscale t_{0}} \prod_{j=1}^{k}
\soms_{j}^{\gpconc-1}.
\end{eqnarray*}
Since the distribution factorizes, the $\{\omsprop_{k}\}$ are
independent of each other and of $\sumofjumps$. Second, we can read off
the distributional kernel of each $\omsprop_{k}$ to establish
$\omsprop
_{k} \stackrel{\mathrm{i.i.d.}}{\sim} \tb(\gpconc,1)$, from whence it follows that
$\stickprop_{k} \stackrel{\mathrm{i.i.d}.}{\sim} \tb(1,\gpconc)$.
\end{pf}
\end{example}
%
%

\subsection{Inference} \label{sec:sub_infer}

In some sense, we skipped ahead in describing inference
in Sections~\ref{sec:epf_infer} and~\ref{sec:stick_infer}.\vadjust{\goodbreak}
There, we made use of the fact that
random labels for partitions and features imply exhangeability
of the data partition block assignments $(\indpart_{n})$ and
data feature assignments $(\indfa_{n})$. In the discussion above,
we study the object that associates random uniformly
distributed labels with each partition or feature. Assuming the
labels come from a uniform distribution rather than a general continuous
distribution is a special case of the discussion in Section~\ref{sec:epf_infer},
and we defer the general case to the next section (Section~\ref{sec:crm}).

We have seen above that it is particularly straightforward
to obtain an EPPF or EFPF formulation, which yields Gibbs sampling
steps as described in Section~\ref{sec:epf_infer}, when the stick lengths
are generated according to a normalized Poisson process in the
partition case or a Poisson process in the feature case. Examples \ref
{ex:sub_ibp}
and~\ref{ex:sub_crp} illustrate how to find such exchangeable probability
functions. Further, we have already seen the usefulness of the stick
representation
in inference, and Examples~\ref{ex:sub_ibp} and~\ref
{ex:sub_crp_sticks} illustrate how
stick-length distributions may be recovered from the subordinator framework.

\section{Completely Random Measures} \label{sec:crm}

In our discussion of subordinators, the jump sizes of the subordinator
corresponded to the feature frequencies or unnormalized partition frequencies
and were the quantities of interest. By contrast, the locations of the jumps
mainly served as convenient labels for the frequencies.
These locations were chosen uniformly at random from the unit interval.
This choice guaranteed the a.s. uniqueness of the labels and the
exchangeability of the sequence of index assignments: $(\indpart_{n})$
in the clustering case or
$(\indfa_{n})$ in the feature case.

However, a labeling retains exchangeability and a.s. uniqueness as long
as the labels are chosen i.i.d. from any continuous distribution (not just the
uniform distribution).
Moreover, in typical applications, we wish to associate some parameter,
often referred to as a ``random effect,'' with each partition block or feature.
In the partition case,
we usually model the $n$th data point $\data_{n}$
as being generated according to some
likelihood depending
on the parameter corresponding to its block assignment.
For example, an individual animal's height and weight, $\data_{n}$,
varies randomly
around the height and weight of its species, $\indpart_{n}$.
Likewise, in the feature
case,
we typically model the observed data point $\data_{n}$
as being generated according to some likelihood
depending on the collection of parameters corresponding
to its collection of feature
block assignments [cf. equation~(\ref{eq:likelihood})]. For example,
the book-buying
pattern of
an online consumer, $\data_{n}$, varies with some noise based on
the topics this person likes to read about: $\indfa_{n}$ is a
collection, possibly empty, of such topics.

In these cases, it can be useful
to suppose that the partition block labels (or feature labels)
$\randlabfa_{k}$ are not
necessarily $\mathbb{R}_{+}$-valued but rather are generated i.i.d.
according to
some continuous distribution $\basedistr$ on a general space $\rmsp$.
Then, whenever $k$ is the order-of-appearance partition block label of index $n$, we let
$\indpart_{n} = \randlabpart_{k}$. Similarly, whenever $k$ is the
order-of-appearance feature label for some feature to which
index $n$ belongs, $\randlabfa_{k} \in\indfa_{n}$. Finally, then,
we complete the generative model in the partition case by letting
$\data_{n} \indep\likedistr(\indpart_{n})$ for some distribution
function $\likedistr$
depending on parameter $\indpart_{n}$. And in the feature case,
$\data_{n} \indep\likedistr(\indfa_{n})$, where now the distribution
function $\likedistr$ depends on the collection of parameters~$\indfa_{n}$.

When we take the jump sizes $(\rmweight_{k})$ of a subordinator as the weights
of atoms with locations $(\rmloc_{k})$ drawn i.i.d. according to
$\basedistr$ as
described above, we find ourselves with a \textit{completely random measure}
$\randm$:
%
%
\begin{equation}
\label{eq:crm} \randm= \sum_{k=1}^{\infty}
\rmweight_{k} \delta_{\rmloc_{k}}.
\end{equation}
A completely random measure is a random measure $\mu$ such that
whenever $A$ and $A'$ are disjoint sets, we have that $\mu(A)$
and $\mu(A')$ are independent random variables.

To see that associating these more general atom locations to the jumps
of a subordinator
yields a completely random measure,
note that Theorem~\ref{thm:sub_ppp} tells us that the subordinator
jump sizes are
generated according to a Poisson point process,
with some intensity measure $\intens(d \rmweight)$.
The Marking Theorem for Poisson point processes \citep
{kingman:1993:poisson} in turn yields that the
tuples $\{(\rmweight_{k}, \rmloc_{k})\}_{k}$ are generated according to
a Poisson point process
with intensity measure $\intens(d \rmweight) \rmlocdist(d\rmloc)$. By
\citet{kingman:1967:completely}, whenever the tuples $\{(\rmweight_{k},
\rmloc_{k})\}_{k}$ are
drawn according to a Poisson point process, the measure in
equation~(\ref{eq:crm}) is
completely random.

%
%
\begin{example}[(Dirichlet process)] \label{ex:crm_dp}
We can form a completely random measure from the
gamma process subordinator and a random labeling of the partition blocks.
Specifically, suppose that the labels come from a continuous measure
$\rmlocdist$. Then
we generate a completely random measure $\gproc$, called a
\textit{gamma process} \citep{Ferguson73},
in the
following way:
%
%
\begin{eqnarray}
\label{eq:gamma_ppp_intens} \nu(d\rmweight\times d
\rmloc) &=& \gpconc\rmweight^{-1} e^{-\gpscale\rmweight} \,d\rmweight\cdot
\rmlocdist(d\rmloc),
\\
\label{eq:gamma_ppp_draw} \bigl\{(
\rmweight_{k},\rmloc_{k})\bigr\}_{k} &\sim&\ppp
(\intens),
\\
\label{eq:gamma_proc_crm} \gproc&=& \sum
_{k=1}^{\infty} \rmweight_{k}
\delta_{\rmloc_{k}}.
\end{eqnarray}
Here, $\ppp(\intens)$ denotes a draw from a Poisson point process with
intensity measure $\intens$. The parameters $\gpconc> 0$ and
$\gpscale
> 0$ are
the same as for the gamma process subordinator. A gamma process draw,
along with its generating Poisson point process intensity measure, is
illustrated in Figure~\ref{fig:gamma_process}.

%
\begin{figure}

\includegraphics{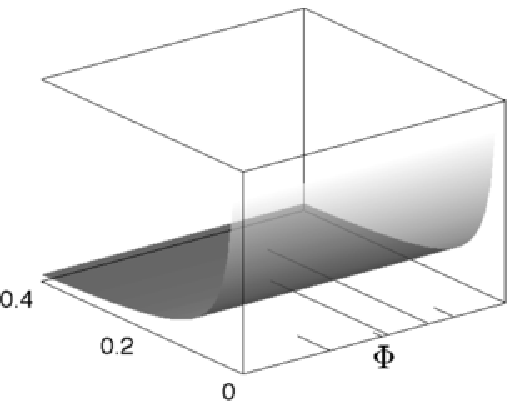}

\caption{The gray manifold depicts the
Poisson point process intensity measure $\intens$ in
equation (\protect\ref{eq:gamma_ppp_intens}) for the choice $\rmsp= [0,1]$
and $\rmlocdist$ the
uniform distribution on $[0,1]$. The endpoints of the line segments are
points drawn from the Poisson point process as in
equation (\protect\ref{eq:gamma_ppp_draw}).
Taking the positive real-valued coordinate (leftmost axis) as the atom
weights, we find the random measure $\gproc$ (a gamma process) on
$\rmsp
$ from equation~(\protect\ref{eq:gamma_proc_crm}) in the bottom
plane.}\label{fig:gamma_process}
\end{figure}

The \textit{Dirichlet process} (DP) is the random measure formed by
normalizing the gamma process \citep{Ferguson73}. Since the Dirichlet
process atom weights
sum to one, it cannot be completely random. We can write the Dirichlet
process $\dproc$
generated from the gamma process $G$ above as
\begin{eqnarray*}
\tau&=& \sum_{k=1}^{\infty}
\rmweight_{k},
\\
\partblockfreq_{k} &= &\rmweight_{k} / \tau,
\\
\dproc&=& \sum_{k=1}^{\infty}
\partblockfreq_{k} \delta _{\randlabpart_{k}}.
\end{eqnarray*}
The random variables $\partblockfreq_{k}$ have the same distribution as
the Dirichlet process sticks
[equation (\ref{eq:dp_stick})]
or normalized gamma process subordinator jump lengths, as we have seen
above (Example~\ref{ex:sub_crp}).
\end{example}

Consider sampling points from a Dirichlet process and
forming the induced partition of the data\vadjust{\goodbreak} indices. Theorem \ref
{thm:eppf_gamma} shows
us that the distribution of the induced partition is the Chinese restaurant
process EPPF.

\begin{figure}

\includegraphics{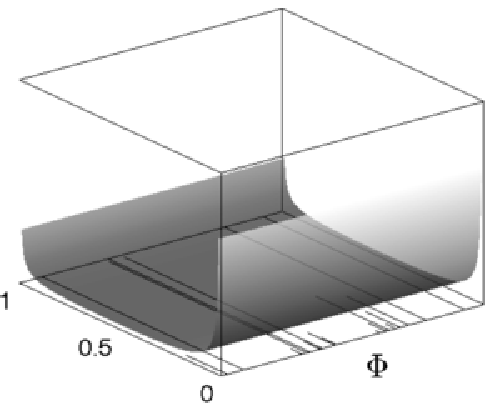}

\caption{The gray manifold depicts the Poisson
point process intensity measure $\intens$ in equation (\protect\ref
{eq:beta_ppp_intens}) for
the choice $\rmsp= [0,1]$ and $\rmlocdist$ the uniform distribution on
$[0,1]$. The endpoints of the line segments are points drawn from the
Poisson point process as in equation (\protect\ref
{eq:beta_ppp_draw}). Taking the
$[0,1]$-valued coordinate (leftmost axis) as the atom weights, we find
the measure $\bproc$ (a~beta process) on $\rmsp$ from
equation (\protect\ref{eq:beta_ppp_crm})
in the bottom plane.}\label{fig:beta_process}
\end{figure}

%
%
\begin{example}[(Beta process)] \label{ex:crm_bp}
We can form a completely random measure from the beta process
subordinator and a random labeling of the feature blocks. If the labels
are generated i.i.d. from a continuous measure $\basedistr$, then we say the
completely random measure $B$, generated as follows, is called a
\textit{beta process}:
%
%
\begin{eqnarray}
\label{eq:beta_ppp_intens} \qquad\nu(d\rmweight\times d
\rmloc) &= &\bpmass\bpconc\rmweight^{-1} (1-\rmweight)^{\bpconc-1} \,d
\rmweight \cdot\rmlocdist(d\rmloc),
\\
\label{eq:beta_ppp_draw} \bigl\{(\rmweight_{k},
\rmloc_{k})\bigr\}_{k} &\sim&\ppp(\intens),
\\
\label{eq:beta_ppp_crm} \bproc&=& \sum
_{k=1}^{\infty} \rmweight_{k}
\delta_{\rmloc_{k}}.
\end{eqnarray}
The beta process, along with its generating intensity measure, is
depicted in Figure~\ref{fig:beta_process}.
The $(\rmweight_{k})$ have the same distribution as the beta process
sticks [equation~(\ref{eq:beta_proc_stick_lengths})]
or the beta process subordinator jump lengths (Example~\ref{ex:sub_ibp}).
\end{example}

%

Now consider sampling a collection of atom locations
according to Bernoulli draws from the atom weights of a beta process and
forming the induced feature allocation of the data indices.
Theorem~\ref{thm:sub_sticks} shows
us that the distribution of the induced feature allocation is given by
the Indian buffet
process EFPF.

\subsection{Inference} \label{sec:crm_infer}

In this section we finally study the full model first outlined in
the context of inference of partition and feature structures in
Section~\ref{sec:epf_infer}. The partition or feature labels\vadjust{\goodbreak} described in
this section are the same as the block-specific parameters first
described
in Section~\ref{sec:epf_infer}. Since this section focuses on a
generalization of the partition or feature labeling scheme beyond
the uniform distribution option encoded in subordinators,
inference for the atom weights remains unchanged
from Sections~\ref{sec:epf_infer},~\ref{sec:stick_infer} and \ref
{sec:sub_infer}.

However, we note that, in the course of inferring
underlying partition or feature structures, we are
often also interested in inferring the parameters of
the generative model of the data given the partition
block or the feature labels. Conditional on the
partition or feature structure, such inference
is handled as in a normal hierarchical model
with fixed dependencies. Namely, the parameter
within a particular block may be inferred from the data points
that depend on this block as well as the prior distribution
for the parameters. Details for the Dirichlet process
example inferred via MCMC sampling are provided by
\citet{maceachern:1994:estimating},
\citet{escobar:1995:bayesian},
\citet{neal:2000:markov};
\citet{blei:2006:variational} work out details for the Dirichlet
process using variational methods.
In the beta process case, \citet{griffiths:2006:infinite},
\citet{teh:2007:stick},
\citet{ThibauxJo07}
describe MCMC sampling, and \citet{paisley:2010:stick}
describe a variational approach.

\section{Conclusion} \label{sec:conclusion}

In the discussion above we have pursued a progressive augmentation from
(1) simple
distributions over partitions and feature allocations in the form of
exchangeable
probability functions to (2) the representation of stick lengths
encoding frequencies
of the partition block and feature occurrences to (3) subordinators,
which associate
random $\mathbb{R}_{+}$-valued labels with each partition block or
feature, and finally
to (4) completely random measures, which associate a general class of
labels with
the stick lengths and whose labels we generally use as parameters in
likelihood models
built from the partition or feature allocation representation.

Along the way, we have focused primarily on two vignettes. We have shown,
via these successive augmentations, that the Chinese restaurant
process specifies the marginal distribution of the induced partition
formed from
i.i.d. draws from a Dirichlet process, which is in turn a normalized
completely random
measure. And we have shown that the Indian buffet process specifies the marginal
distribution of the induced feature allocation formed by i.i.d. Ber\-noulli
draws across
the weights of a beta process.

There are many extensions of these ideas that lie beyond the scope of
this paper.
A number of extensions of the CRP and Dirichlet process exist---in
either the EPPF form (\cite{pitman:1996:some}; \cite{blei:2010:distance}),
the stick-length form (Dunson and Park, \citeyear{dunson:2008:kernel})
or the random measure form (Pitman and Yor, \citeyear{pitman:1997:two}). Likewise,
extensions of
the IBP and beta process have been explored (\cite
{teh:2007:stick}; \cite{paisley:2010:stick}; \cite{broderick:2012:beta}).

More generally, the framework above demonstrates how alternative partition
and feature allocation models may be constructed---either by introducing
different EPPFs (\cite{pitman:1996:some}; \cite{gnedin:2006:exchangeable}) or EFPFs,
different stick-length distributions \citep{ishwaran:2001:gibbs}
or different random measures \citep{wolpert:2004:reflecting}.

Finally, we note that expanding the set of combinatorial structures
with useful
Bayesian priors from partitions to the superset of feature allocations suggests
that further such structures might be usefully \mbox{examined}. For instance, the
\textit{beta negative binomial process} (\cite
{broderick:2011:combinatorial}; \cite{zhou:2012:beta})
provides a prior on a generalization of a
feature allocation where we allow the features themselves to be
multisets; that is,
each index may have nonnegative integer multiplicities of features. Models
on trees (\cite{adams:2010:tree}; \cite{mccullagh:2008:gibbs};
\cite{blei:2010:nested}),
graphs \citep{li:2006:pachinko} and permutations \citep
{pitman:1996:some} provide avenues for future exploration.
And there likely remain further structures to be fitted out with
useful Bayesian priors.


\section*{Acknowledgments}
T. Broderick's research was funded by a National Science Foundation
Graduate Research Fellowship.
This material is supported in part by the National Science Foundation
Award 0806118 Combinatorial
Stochastic Processes and is based upon work supported in part by the
Office of
Naval Research under contract/grant number N00014-11-1-0688.

%

%

\end{document}